\documentclass[12pt, a4paper]{amsart}
\usepackage{fullpage}
\usepackage{amscd,amsmath,amsthm,amssymb,graphics,color}
\usepackage[colorlinks,linkcolor=blue,filecolor=blue,citecolor=blue,urlcolor=blue,bookmarks=false,hypertexnames]{hyperref} 
\usepackage{tikz}
\usepackage{caption,subcaption}

\renewcommand{\Bbb}{\mathbb}

\def\NZQ{\Bbb}               
\def\NN{{\NZQ N}}

\def\ZZ{{\NZQ Z}}

\def\FF{{\NZQ F}}

\def\frk{\mathfrak}               

\def\mm{{\frk m}}

\def\Phi{{\frk n}}
\def\Phi{{\frk N}}

\def\MI{{\mathcal I}}

\def\ab{{\bold a}}

\def\opn#1#2{\def#1{\operatorname{#2}}} 
\opn\chara{char} \opn\length{\ell} \opn\pd{pd} \opn\rk{rk}
\opn\projdim{proj\,dim} \opn\injdim{inj\,dim} \opn\rank{rank}
\opn\depth{depth} \opn\grade{grade} \opn\height{height}\opn\coheight{coheight}
\opn\embdim{emb\,dim} \opn\codim{codim}

\opn\Tr{Tr} \opn\bigrank{big\,rank}
\opn\superheight{superheight}\opn\lcm{lcm}
\opn\trdeg{tr\,deg}
\opn\reg{reg} \opn\lreg{lreg} \opn\ini{in} \opn\lpd{lpd}
\opn\size{size}\opn\bigsize{bigsize}
\opn\cosize{cosize}\opn\bigcosize{bigcosize}
\opn\sdepth{sdepth}\opn\sreg{sreg}
\opn\link{link}\opn\fdepth{fdepth}\opn\type{type}
\opn\GL{GL} 
 \opn\width{width}

\opn\div{div} \opn\Div{Div} \opn\cl{cl} \opn\Cl{Cl}
\opn\Spec{Spec} \opn\Supp{Supp} \opn\supp{supp} \opn\Sing{Sing}
\opn\Ass{Ass} \opn\Min{Min}\opn\Mon{Mon} \opn\dstab{dstab} \opn\astab{astab}
\opn\Syz{Syz}
\opn\Ann{Ann} \opn\Rad{Rad} \opn\Soc{Soc}
\opn\Im{Im} \opn\ker{ker} \opn\Coker{Coker} \opn\Am{Am}
\opn\Hom{Hom} \opn\Tor{Tor} \opn\Ext{Ext} \opn\End{End}
\opn\Aut{Aut} \opn\id{id}

\opn\nat{nat}
\opn\pff{pf}
\opn\Pf{Pf} \opn\GL{GL} \opn\SL{SL} \opn\mod{mod} \opn\ord{ord}
\opn\Gin{Gin} \opn\Hilb{Hilb}\opn\sort{sort}
\opn\Proj{Proj} 
\opn\aff{aff} \opn\con{conv} \opn\relint{relint} \opn\st{st}
\opn\lk{lk} \opn\cn{cn} \opn\core{core} \opn\vol{vol}
\opn\link{link} \opn\star{star}\opn\lex{lex}
\opn\gr{gr}

\opn\dirlim{\underrightarrow{\lim}}
\opn\inivlim{\underleftarrow{\lim}}
\let\dirsum=\oplus
\let\tensor=\otimes
\let\iso=\cong

\let\Dirsum=\bigoplus

\let\ab\allowbreak
\let\to=\rightarrow

\def\Implies{\ifmmode\Longrightarrow \else
	\unskip${}\Longrightarrow{}$\ignorespaces\fi}
\def\implies{\ifmmode\Rightarrow \else
	\unskip${}\Rightarrow{}$\ignorespaces\fi}
\def\iff{\ifmmode\Longleftrightarrow \else
	\unskip${}\Longleftrightarrow{}$\ignorespaces\fi}

\let\:=\colon
\newtheorem{Theorem}{Theorem}[section]
\newtheorem{Lemma}[Theorem]{Lemma}
\newtheorem{Corollary}[Theorem]{Corollary}
\newtheorem{Proposition}[Theorem]{Proposition}
\newtheorem{Remark}[Theorem]{Remark}

\newtheorem{Example}[Theorem]{Example}

\newtheorem{Definition}[Theorem]{Definition}

\newtheorem{Notation}[Theorem]{Notation}

\newtheorem{Discussion}[Theorem]{Discussion}

\newenvironment{discussionbox}[1][]{%
	\begin{Discussion}[#1]\pushQED{\qed}}{\popQED \end{Discussion}}
\newenvironment{notationbox}[1][]{%
	\begin{Notation}[#1]\pushQED{\qed}}{\popQED \end{Notation}}

\let\epsilon\varepsilon

\title {On the linearity of the syzygies of Hibi rings}

\author {Dharm Veer}

\dedicatory{In memory of Prof. C.~S.~Seshadri}
\address{Chennai Mathematical Institute, Siruseri, Tamilnadu 603103, India.}
\email{dharm@cmi.ac.in}

 \thanks{The author was partly supported by the grant CRG/2018/001592 (Manoj Kummini) from Science and Engineering Research Board, India and by an Infosys Foundation fellowship.}


\begin{document}

\begin{abstract}
In this article, we prove necessary conditions for Hibi rings to satisfy Green-Lazarsfeld property $N_p$ for $p=2$ and $3$. We also show that if a Hibi ring satisfies property $N_4$, then it is a polynomial ring or it has a linear resolution. Therefore, it satisfies property $N_p$ for all $p\geq 4$ as well. As a consequence, we characterize distributive lattices whose comparability graph is chordal in terms of the subposet of join-irreducibles of the distributive lattice. Moreover, we characterize complete intersection Hibi rings. 
\end{abstract} 

\maketitle

\section{Introduction}

A classical problem in commutative algebra is to study graded minimal free resolutions of graded modules over polynomial rings. 
Let $S$ be a standard graded polynomial ring in finitely many variables over a field $K$ and $I$ be a graded $S$-ideal generated by quadratics. 
To study the graded minimal free resolution of $S/I$, Green-Lazarsfeld \cite{[GL86]} defined property $N_p$ for $p\in \NN$. The ring $S/I$ satisfies property $N_p$ if $S/I$ is normal and the graded minimal free resolution of $S/I$ over $S$ is linear upto $p^{th}$ position. In this article, we study Green-Lazarsfeld property $N_p$ for Hibi rings. We also characterize complete intersection Hibi rings. 

Let $L$ be a finite distributive lattice and $P=\{p_1,\ldots,p_{n}\}$ be the subposet of join-irreducible elements of $L$. Let $K$ be a field and let $R=K[t,z_1,\ldots,z_n]$ be a polynomial ring over $K$. The Hibi ring associated with $L$, denoted by $R[L]$,  is the subring of $R$ generated by the  monomials $u_{\alpha}=t\prod_{p_i\in \alpha}z_i$ where $\alpha\in L$. Hibi \cite{[HIBI87]} showed that $R[L]$ is a normal Cohen–Macaulay
domain of dimension $\#P+1$, where $\#P$ is the cardinality of $P$. See~\cite{[ENE15]} for a survey on Hibi rings. If we set $\deg(t)=1$ and $\deg(z_i)=0$ for all $1\leq i \leq n$, then $R[L]$ may be viewed as a standard graded algebra over $K$. Let $K[L]=K[\{x_\alpha: \alpha\in L\}]$ be the polynomial ring over $K$ and $\pi: K[L] \to R[L]$ be the $K$-algebra homomorphism with $x_\alpha \mapsto u_{\alpha}$. The ideal $I_L=(x_\alpha x_\beta-x_{\alpha\wedge\beta}x_{\alpha\vee\beta}: \alpha,\beta\in L\ \text{and}\ \alpha,\beta \text{ incomparable})$ is the kernel of the map $\pi$. It is called the Hibi ideal associated to $L$.\par

In past, various authors have studied minimal free resolution of Hibi rings \cite{[DM17],[EHH15],[EHM15],[ENE15]}. In Section~\ref{propertyn2} and \ref{propertn3}, we establish some results about property $N_p$ of Hibi rings [Theorem~\ref{planar}, \ref{connected}, \ref{beta35connected} and \ref{propertyn4}]. In \cite{[EHH15]}, the authors have classified all Hibi rings which have linear resolution. In \cite{[EHM15]}, the authors have given a combinatorial description of regularity of Hibi rings in terms of poset of join-irreducibles. Ene \cite{[ENE15]} characterizes all simple planar distributive lattices for which the associated Hibi ring satisfies property $N_2$ [Theorem~\ref{ene}]. The Segre product of polynomial rings may be viewed as a Hibi ring. Property $N_p$ for Segre product of polynomial rings has been completely studied for all $p$, see \cite{VEE22} for chronological developments. In~\cite{VEE22}, the author has proved that if a Hibi ring satisfies property $N_2$, then its Segre product with a polynomial ring in finitely many variables also satisfies property $N_2$. When the polynomial ring is in two variables, the above statement was proved for $N_3$. The author has also studied minimal Koszul syzygies of Hibi ideals and initial Hibi ideals in \cite{VEE22}.

Let $P$ be a poset. The comparability graph $G_P$ of $P$ is a graph on the underlying set of $P$ such that $\{x,y\}$ is an edge of $G_P$ if and only if $x$ and $y$ are comparable in $P$.
Hibi and Ohsugi~\cite{OH17chordalcomparabilitygraphs} characterized chordal comparability graph of posets using toric ideals associated with multichains of poset. Using one of our results and~\cite[Theorem\ 1]{Fro90}, we characterize chordal comparability graph of distributive lattices in terms of the subposet of join-irreducibles of the distributive lattice~[Corollary~\ref{coro:chordalcomparabilitygraph}].

The article is organised as follows. In Section~\ref{hibiringdefi}, we recall some basic notions of algebra and combinatorics. In~\cite{[EHH15]}, the authors have introduced the notion of homologically pure subsemigroup of an affine semigroup, and proved that if $H'$ is a homologically pure subsemigroup of an affine semigroup $H$, then the Betti numbers of the semigroup ring $K[H]$ are greater than equal to those of $K[H']$~[Proposition~\ref{homologicallypure}]. In Section~\ref{sec:homologicallypure}, we identify two kinds of homologically pure subsemigroups of an affine semigroup associated to a Hibi ring. Using these, we prove some sufficient conditions for Hibi rings to not satisfy property $N_2$ in Section~\ref{propertyn2}~[Theorem~\ref{planar}, \ref{connected}].\par 

In Section~\ref{propertn3}, we study property $N_p$ of Hibi rings for $p\geq 3$. First, we prove that if a poset is connected and it has atleast two minimal and atleast two maximal elements,  then the associated Hibi ring does not satisfy property $N_3$ [Theorem~\ref{beta35connected}]. The second main result of this section is that a Hibi ring satisfies property $N_4$ if and only if either it is a polynomial ring or it has a linear resolution if and only if it is a polynomial ring or its initial ideal has a linear resolution~[Theorem~\ref{propertyn4}]. We also characterize all such Hibi rings combinatorially which gives a different proof of~\cite[Corollary\ 10]{EQR13regularityplanardistlattices}. 
The last section of the article is devoted to combinatorial characterization of complete intersection Hibi rings~[Theorem~\ref{completeinter}].


\subsection*{Acknowledgements}
I am extremely grateful to Manoj Kummini for his guidance and various insightful discussions throughout the preparation of this article. The computer algebra systems Macaulay2~\cite{[M2]} and SageMath~\cite{[sagemath]} provided valuable assistance in studying examples.

\section{Preliminaries}\label{hibiringdefi}
 
We start by defining some notions of posets and distributive lattices. For more details and examples, we refer the reader to \cite[Chapter\ 3]{[STAN12]}. Throughout this article, all posets and distributive lattices will be finite.
 
Let $P$ be a poset. We say that two elements $x$ and $y$ of $P$ are {\em comparable}  if $x\leq y$ or $y\leq x$; otherwise $x$ and $y$ are {\em incomparable}. For $x,y \in P$, we say that {\em y covers x} if $x<y$ and there is no $z\in P$ with $x<z<y$. We denote it by $x\lessdot y$. A poset is completely determined by its cover relations. The {\em Hasse diagram} of poset $P$ is the graph whose vertices are elements of $P$, whose edges are cover relations, and such that if $x<y$ then $y$ is ``above'' $x$ (i.e. with a higher vertical coordinate). In this article, we use the Hasse diagrams to represent posets. 

A {\em chain} $C$ of $P$ is a totally ordered subset of $P$. The {\em length} of a chain $C$ of $P$ is $\#C - 1$. The {\em rank} of $P$, denoted by $\rank(P)$, is the maximum of the lengths of chains in $P$.
A poset is called {\em pure} if its all maximal chains have the same length. A subset $\alpha$ of $P$ is called an {\em order ideal} of $P$ if it satisfies the following condition: for any $x \in \alpha$ and $y \in P$, if $y\leq x$, then  $y\in \alpha$. Define $\MI(P):= \{\alpha \subseteq P: \alpha  \ \text{is an order ideal of}\ P \}$. It is easy to see that $\MI(P)$, ordered by inclusion, is a distributive lattice under union and intersection. $\MI(P)$ is called the {\em ideal  lattice} of the poset $P$. Let $L$ be a lattice. An element $x\in L$ is called {\em join-irreducible} if $x$ is not the minimal element of $L$ and whenever $x =y \vee z$ for some $y,z \in L$, we have either $x=y$ or $x=z$.

 \begin{Definition}\normalfont
 	Let $P$ and $Q$ be two posets.
 	\begin{itemize}
 		\item[$(a)$] A nonempty subset $S$ of $P$ is an {\em antichain} in $P$ if any two distinct elements of $S$ are incomparable. An antichain with $n$ elements is said to have {\em width} $n$. Define $\width(P):= max\{\#S:S\subseteq P, S\ \text{is an antichain in}\ P\}$.
 		\item[$(b)$]A poset $P$ is called {\em simple} if there is
 		no $p \in P$ with the property that all elements of $P$ are comparable to $p$.
 		\item[$(c)$] The {\em ordinal sum} $P\dirsum Q$ of the disjoint posets $P$ and $Q$ is the poset on the set $P\cup Q$ with the following order: if $x, y \in P\dirsum Q $, then $x \leq y$ if either $x, y \in P$ and $x \leq y$ in $P$ or $x, y \in Q$ and $x \leq y$ in $Q$ or $x \in P$ and $y \in Q$.
 		\item[$(d)$] Let $P,Q$ be two posets on disjoint sets. The {\em disjoint union} of posets $P$ and $Q$ is the poset $P+Q$ on the set $P\cup Q$ with the following order: if $x,y \in P+Q$, then $x\leq y$ if either $x, y \in P$ and $x \leq y$ in $P$ or $x, y \in Q$ and $x \leq y$ in $Q$. A poset which can be written as disjoint union of two posets is called {\em disconnected}. Otherwise, $P$ is {\em connected}.
 		\item[$(e)$] $P$ and $Q$ are said to be {\em isomorphic}, denoted by $P \cong Q$, if there exists an order-preserving bijection $\varphi: P \to Q$  whose inverse is order preserving.		
 		
 		\item[$(f)$] A subposet $P'$ of $P$ is said to be a {\em cover-preserving subposet} of $P$ if for every $x,y \in P'$ with $x\lessdot y$ in $P'$, we have $x\lessdot y$ in $P$.
 	\end{itemize}
 \end{Definition}
 
 \begin{Example}
 	Let $P$ be the poset as shown in Figure~\ref{fig1}. Let $P'$ and $P''$ be the subposets of $P$ as shown in Figure~\ref{fig2} and Figure~\ref{fig3} respectively. It is easy to see that $P'$ is a cover-preserving subposet of $P$ but $P''$ is not  a cover-preserving subposet of $P$ since $p_3\lessdot p_7$ in $P''$ but not in $P$.
 \end{Example}

 \begin{figure}[h]
 	\begin{subfigure}[t]{4cm}
 		\centering	
 		\begin{tikzpicture}[scale=2]
 		\draw[] (0,.5)--(-.5,0)--(-.5,.5)--(0,1) (0,0)--(0,.5)--(0,1)--(0,1.5) (0,1)--(.5,1.5) (1,1.5)--(1,1)--(1,.5)--(.5,0)--(.5,.5)--(.5,1)--(.5,1.5) (0,.5)--(1,1);

 		\filldraw[black] (0,0) circle (.5pt) node[anchor=north]  {$p_1$}; 
 		\filldraw[black] (.5,0) circle (.5pt) node[anchor=north] {$p_2$}; 
 		\filldraw[black] (0,.5) circle (.5pt) node[anchor=east] {$p_3$}; 
 		\filldraw[black] (.5,.5) circle (.5pt) node[anchor=west] {$p_4$}; 
 		\filldraw[black] (0,1) circle (.5pt) node[anchor=east]  {$p_5$}; 
 		\filldraw[black] (.5,1) circle (.5pt) node[anchor=west] {$p_6$}; 
 		\filldraw[black] (0,1.5) circle (.5pt) node[anchor=east] {$p_7$}; 
 		\filldraw[black] (.5,1.5) circle (.5pt) node[anchor=west] {$p_8$}; 
 		\filldraw[black] (-.5,.5) circle (.5pt) node[anchor=east] {$p_9$}; 
 		\filldraw[black] (-.5,0) circle (.5pt) node[anchor=east] {$p_{10}$}; 
 		\filldraw[black] (1,.5) circle (.5pt) node[anchor=west]  {$p_{11}$}; 
 		\filldraw[black] (1,1) circle (.5pt) node[anchor=west] {$p_{12}$}; 
 		\filldraw[black] (1,1.5) circle (.5pt) node[anchor=west] {$p_{13}$}; 
 		\end{tikzpicture}
 		\caption{}\label{fig1}
 	\end{subfigure}
 	\quad
 	\begin{subfigure}[t]{4cm}
 		\centering	
 		\begin{tikzpicture}[scale=2]
 		\draw[] (0,0)--(0,.5)--(0,1)--(0,1.5) (.5,0)--(.5,.5)--(.5,1)--(.5,1.5) (.0,1)--(.5,1.5);

 		\filldraw[black] (0,0) circle (.5pt) node[anchor=north]  {$p_1$}; 
 		\filldraw[black] (.5,0) circle (.5pt) node[anchor=north] {$p_2$}; 
 		\filldraw[black] (0,.5) circle (.5pt) node[anchor=east] {$p_3$}; 
 		\filldraw[black] (.5,.5) circle (.5pt) node[anchor=west] {$p_4$}; 
 		\filldraw[black] (0,1) circle (.5pt) node[anchor=east]  {$p_5$}; 
 		\filldraw[black] (.5,1) circle (.5pt) node[anchor=west] {$p_6$}; 
 		\filldraw[black] (0,1.5) circle (.5pt) node[anchor=south] {$p_7$}; 
 		\filldraw[black] (.5,1.5) circle (.5pt) node[anchor=south] {$p_8$}; 
 		\end{tikzpicture}
 		\caption{}\label{fig2}
 	\end{subfigure}
 	\quad
 	\begin{subfigure}[t]{4cm}
 		\centering	
 		\begin{tikzpicture}[scale=2]
 		\draw[] (.5,.5)--(0,0)--(0,.5)--(0,1) (.5,0)--(.5,.5)--(.5,1) (.5,.5)--(0,1) (.5,1)--(0,.5);

 		\filldraw[black] (0,0) circle (.5pt) node[anchor=north]  {$p_{10}$}; 
 		\filldraw[black] (.5,0) circle (.5pt) node[anchor=north] {$p_1$}; 
 		\filldraw[black] (0,.5) circle (.5pt) node[anchor=east] {$p_9$}; 
 		\filldraw[black] (.5,.5) circle (.5pt) node[anchor=west] {$p_3$}; 
 		\filldraw[black] (0,1) circle (.5pt) node[anchor=south]  {$p_7$}; 
 		\filldraw[black] (.5,1) circle (.5pt) node[anchor=south] {$p_8$}; 
 		
 		\end{tikzpicture}
 		\caption{}\label{fig3}
 	\end{subfigure}	
 	\caption{}
 \end{figure}
 
\subsection{Hibi rings}\label{subsec:hibirings} Let $L$ be a distributive lattice and let $P$ be the subposet of join-irreducible elements of $L$. By Birkhoff's fundamental structure theorem \cite[Theorem\ 3.4.1]{[STAN12]}, $L$ is isomorphic to the ideal lattice $\MI(P)$. Write $P=\{p_1,\ldots,p_n\}$ and let $R=K[t,z_1,\ldots,z_n]$ be a polynomial ring in $n+1$ variables over a field $K.$ The {\em Hibi ring} associated with $L$, denoted by $R[L]$,  is the subring of $R$ generated by the  monomials $u_{\alpha}=t\prod_{p_i\in \alpha}z_i$ where $\alpha\in L$.  If we set $\deg(t)=1$ and $\deg(z_i)=0$ for all $1\leq i \leq n$, then $R[L]$ may be viewed as a standard graded algebra over $K$.

Let $K[L]=K[\{x_\alpha: \alpha\in L\}]$ be the polynomial ring over $K$ and $\pi: K[L] \to R[L]$ be the $K$-algebra homomorphism with $x_\alpha \mapsto u_{\alpha}$. Let \[I_L=(x_\alpha x_\beta-x_{\alpha\wedge\beta}x_{\alpha\vee\beta}: \alpha,\beta\in L\ \text{and}\ \alpha,\beta \text{ incomparable})\] be an $K[L]$-ideal. Let $<$ be a total order on the variables of $K[L]$ with the property that one has $x_\alpha < x_\beta$ if $\alpha < \beta$ in $L$. Consider the graded reverse lexicographic order $<$ on $K[L]$ induced by this order of the variables. The generators of $I_L$ described above forms a Gr\"obner basis of $\ker(\pi)$ with respect to $<$~\cite[Theorem\ 6.19]{[HHO18]}. In particular, $\ker(\pi) = I_L$.
The ideal $I_L$ is called the {\em Hibi ideal} of $L$. The {\em initial Hibi  ideal} is
\[\ini_<(I_L)=(x_\alpha x_\beta: \alpha,\beta\in L \ \text{and}\ \alpha, \beta \text{ incomparable}).\]

We now discuss how Hibi rings behave under the ordinal sum of posets. An another result about ordinal sum of two posets will be proved in Lemma~\ref{p_1}.

\begin{Lemma}\label{ordinalsum}
	Let $P_1, \{p\}$ and $P_2$ be posets. Let $P$ be the ordinal sum $P_1\dirsum \{p\}\dirsum P_2$. Then
	$$R[\MI(P)]  \cong  R[\MI(P_1\dirsum P_2)]\tensor_{K}K[y] \cong  R[\MI(P_1)] \tensor_{K}R[\MI(P_2)],$$
	where $K[y]$ is a polynomial ring.
\end{Lemma}

\begin{proof}
	First, we  prove that $$R[\MI(P)]\cong R[\MI(P_1\dirsum P_2)] \tensor_{K}K[y].$$
	Let $ R[\MI(P_1\dirsum P_2)]= K[\{u_\beta\ : \beta \in \MI(P_1\dirsum P_2) \}]/I_{\MI(P_1\dirsum P_2)}$ and $R[\MI(P)]= K[\{v_\alpha\ : \alpha \in \MI(P) \}]/I_{\MI(P)}$.  Define a map
	$$\varphi : K[\{v_\alpha\ : \alpha \in \MI(P) \}]\to T := K[y,\{u_\beta\ : \beta \in \MI(P_1\dirsum P_2) \}]$$ by  
	\[ \varphi(v_\gamma) = \begin{cases} 
	u_\gamma & \text{if} \quad \gamma \subseteq P_1,  \\
	y & \text{if} \quad \gamma = P_1\cup \{p\}, \\
	u_{\gamma'} & \text{if} \quad \gamma = P_1 \cup \{p\}\cup \gamma',\ \text{where} \
	\gamma' \subseteq P_2 . \\
	\end{cases}
	\] 
	It is easy to see that $\varphi$ is an isomorphism. If $\alpha, \beta \in \MI(P)$ are incomparable, then either $\alpha, \beta \in \MI(P_1)$ or $\alpha = P_1 \cup \{p\} \cup \alpha'$ and $\beta = P_1\cup \{p\} \cup \beta'$ where $\alpha',\beta' \in \MI(P_2)$ and $\alpha',\beta'$ incomparable. Let $T' = T/(I_{\MI(P_1\dirsum P_2)}T)$ and $\pi : T \to T'$ be the natural surjection. Thus, $\pi \circ \varphi :  K[\MI(P)]\to T'$ and $\ker(\pi \circ \varphi) = \varphi^{-1}I_{\MI(P_1\dirsum P_2)}T$.\par
	
	It is sufficient to show that  $\varphi (I_{\MI(P)}) = I_{\MI(P_1\dirsum P_2)}T$. Let $\alpha, \beta$ be two incomparable elements of $\MI(P)$. If $\alpha, \beta \in \MI(P_1)$ then $\varphi(v_\alpha v_\beta-v_{\alpha\cap\beta}v_{\alpha\cup\beta})= u_\alpha u_\beta-u_{\alpha\cap\beta}u_{\alpha\cup\beta} \in I_{\MI(P_1)}T$. If $\alpha = P_1 \cup \{p\} \cup \alpha'$ and $\beta = P_1 \cup \{p\} \cup \beta'$ where $\alpha',\beta' \in \MI(P_2)$, then $\varphi(v_\alpha v_\beta-v_{\alpha\cap\beta}v_{\alpha\cup\beta})= u_{\alpha'} u_{\beta'}-u_{{\alpha'}\cap{\beta'}}u_{{\alpha'}\cup{ \beta'}} \in I_{\MI(P_2)}T$. Hence, $\varphi (I_{\MI(P)}) \subseteq I_{\MI(P_1)}T  + I_{\MI(P_2)}T$. 
	
	On the other hand, if $\alpha, \beta$ are two incomparable elements of $\MI(P_1)$ then $\varphi(v_\alpha v_\beta-v_{\alpha\cap\beta}v_{\alpha\cup\beta})= u_\alpha u_\beta-u_{\alpha\cap\beta}u_{\alpha\cup\beta}$ while if $\alpha', \beta'$ are two incomparable elements of $\MI(P_2)$ then $\varphi(v_{P_1 \cup \{p\} \cup \alpha'} v_{P_1\cup\{p\} \cup\beta'}- v_{({P_1\cup\{p\} \cup\alpha'})\cap({P_1\cup\{p\} \cup\beta'})} v_{({P_1\cup\{p\} \cup\alpha'})\cup({P_1\cup\{p\} \cup\beta'})})= u_{\alpha'} u_{\beta'}-u_{\alpha'\cap\beta'}u_{\alpha'\cup\beta'}$. Hence the equality.
	
	
	The minimal generating set of the Hibi ideal $I_{\MI(P)}$ can be partitioned between two disjoint set of variables $\{v_\alpha\ : \alpha \in \MI(P)\ \text{and}\ \alpha \subseteq P_1\}$ and $\{v_\alpha\ : \alpha \in \MI(P)\ \text{and}\ P_1\cup\{p\} \subseteq \alpha\}$. So the Hibi ring $R[\MI(P)]$ admits a tensor product decomposition,  where one of the rings is isomorphic to $R[\MI(P_1)]$ and the other ring is isomorphic to $R[\MI(P_2)]$.
\end{proof}

In \cite{[HIBI87]}, Hibi proved that $R[\MI(P_1)\dirsum \MI(P_2)]  \cong  R[\MI(P_1)] \tensor_{K}R[\MI(P_2)]$. One can immediately check that the poset of join-irreducibles of $\MI(P_1)\dirsum \MI(P_2)$ is isomorphic to $P_1\dirsum \{p\}\dirsum P_2$. 

\begin{Corollary}\label{nonsimpleelements}
	Let $P$ be a poset and $P'= \{p_{i_1},...,p_{i_r}\}$ be the subset of all elements of $P$ which are comparable to every element of $P$. Let $P''$ be the induced subposet of $P$ on the set  $P\setminus P'$. Then,
	$$R[\MI(P)] \cong R[\MI(P'')] \tensor_{K}K[y_1,\ldots,y_r],$$
	where $K[y_1,\ldots,y_r]$ is a polynomial ring.
\end{Corollary}
\begin{proof}
	Without loss of generality, we may assume that $p_{i_1}<\cdots<p_{i_r}$ in $P$. Let $P_0 =\{p\in P : p<p_{i_1}\}$, $P_j=\{p\in P : p_{i_j}<p< p_{i_{j+1}}\}$ for $1<j<r-1$ and $P_r=\{p\in P : p>p_{i_r}\}$. Then $P$ is the ordinal sum $P_0\dirsum \{p_{i_1}\}\dirsum P_1\dirsum\cdots\dirsum \{p_{i_r}\}\dirsum P_r$. Now, the result follows from Lemma~\ref{ordinalsum}.
\end{proof}


\subsection{Green-Lazarsfeld property}\label{propertynp} 
 Let $S=K[x_1,\ldots, x_n]$ be a standard graded polynomial ring over a field $K$ and $I$ be a graded $S$-ideal. Let $\FF$ be the graded minimal free resolution of $S/I$ over $S$:
$$\FF : 0 \to \mathop{\Dirsum}_{j}^{} S(-j)^{\beta_{rj}}\to \cdots \to \mathop{\Dirsum}_{j}^{} S(-j)^{\beta_{1j}} \to \mathop{\Dirsum}_{j}^{} S(-j)^{\beta_{0j}}.$$
The numbers $\beta_{ij}$ are called the minimal graded Betti numbers of $S/I$ over $S$.\par

Let $p\in \NN$. Under the notations as above, we say that $S/I$ satisfies {\em Green-Lazarsfeld property $N_p$} if $S/I$ is normal and $\beta_{ij}(S/I)=0$ for all $i \neq j+1$ and $1\leq i\leq  p$. Therefore, $S/I$ satisfies property $N_0$ if and only if it is normal; it satisfies property $N_1$ if and only if it is normal  and $I$ is generated by quadratics; it satisfies property $N_2$ if and only if it satisfies property $N_1$ and $I$ is linearly presented and so on. 

We know that Hibi rings are normal and Hibi ideals are generated by quadratics. Hence, Hibi rings satisfy property $N_1$. Hibi rings are algebras with straightening laws (ASL) and straightening relations are quadratic \cite[\S \ 2]{[HIBI87]}. ASL with quadratic straightening relations are Koszul \cite{[KEMPH90]}. So for a poset $P$, $\beta_{2j}(R[\MI(P)])= 0$ for  all $j\geq 5$ by \cite[Lemma\ 4]{[KEMPH90]}. Therefore, $R[\MI(P)]$ satisfies property $N_2$ if and only if $\beta_{24}(R[\MI(P)])= 0$. Also, it follows from \cite[Theorem\ 6.1]{[ACI15]} that if $R[\MI(P)]$ satisfies property $N_2$ and $\beta_{35}(R[\MI(P)])= 0$, then it satisfies property $N_3$

\subsection{Planar distributive lattices}\label{planarsection}

In this subsection, we define the notion of planar distributive lattice. We state a result of Ene which characterizes all simple planar distributive lattices for which the associated Hibi ring satisfies property $N_2$.
 
\begin{Definition}\cite[Section\ 6.4]{[HHO18]}
A finite distributive lattice $L = \MI(P)$ is called {\em planar} if P can be decomposed into a disjoint union $P = \{p_1,\ldots,p_m \}\cup \{q_1,\ldots,q_n\}$, where $m,n \geq 0$ such that $\{p_1,\ldots,p_m \}$ and $\{q_1,\ldots,q_n \}$ are chains in $P$.
\end{Definition}

\begin{Remark}\normalfont
Let us consider the infinite distributive lattice $\NN^2$ with the
partial order defined as $(i, j) \leq (k, l)$ if $i \leq k$ and $j \leq l$. Let $L = \MI(P)$ be a finite planar distributive lattice, where $P = \{p_1,\ldots,p_m \}\cup \{q_1,\ldots,q_n\}$. Assume that $\{p_1,\ldots,p_m \}$ and $\{q_1,\ldots,q_n \}$ are chains in $P$ with $p_1 \leq \cdots \leq p_m$ and $q_1\leq \cdots \leq q_n$. Define a map 	$$\varphi: \MI(P) \to \NN^2$$
 by 
 	\[ \varphi(\alpha) = \begin{cases} 
 	(0,0) & \text{if} \quad \alpha = \emptyset\text{,}  \\
 	(i,0) & \text{if} \quad \alpha = \{p\in P :\ p\leq p_i\}\text{,} \\
 	(0,j) & \text{if} \quad \alpha = \{p\in P :\ p\leq q_j\}\ \text{,} \\
 	(i,j) & \text{if} \quad \alpha = \{p\in P :\text{either}\ p\leq p_i \ \text{or} \ p\leq q_j \}\text{.} \\
 	\end{cases}
 	\]
 It is easy to see that $\varphi$ is an order-preserving injective map. Hence, any finite planar distributive lattice can be embedded into $\NN^2$. Also, observe that  $[(0, 0), (m,n)]$ is the smallest interval of $\NN^2$ which contains $L$.
\end{Remark}

Let $L$ be a distributive lattice. If the poset of join-irreducibles of $L$ is a simple poset, then sometimes we abuse the notation and say that $L$ is a {\em simple distributive lattice}.

\begin{Theorem}\cite[Theorem\ 3.12]{[ENE15]}\label{ene}
	Let $L= \MI(P)$ be a simple planar distributive lattice with $\# P = n + m$, $L\subset [(0,0),(m,n)]$ with $m,n\geq 2$. Then $R[\MI(P)]$ satisfies property $N_2$ if and only if  the following conditions hold:
	\begin{itemize}
		\item [(i)] At least one of the vertices $(m,0)$ and $(0,n)$ belongs to $L$.
		\item [(ii)] The vertices $(1,n-1)$ and $(m-1,1)$ belong to $L.$
	\end{itemize}
\end{Theorem}

\begin{Corollary}\label{N2planarposets}
	Let $L=\MI(P)$ be a simple planar distributive lattice with $P= \{a_1,\ldots,a_m,\ab b_1,\ldots,b_n\}$. Let  $\{a_1,\ldots,a_m\}$ and $\{b_1,\ldots,b_n\}$ be chains in $P$ with $a_1\lessdot a_2\lessdot \cdots\lessdot a_m$ and $b_1\lessdot b_2\lessdot \cdots\lessdot b_n$. Assume that $\{a_1,\ldots,a_m\}$ is an order ideal of $P$. If $R[\MI(P)]$ satisfies property $N_2$, then $P$ is one of the posets as shown in Figure~\ref{fig20}.
\end{Corollary}

\begin{figure}[h]
	\begin{subfigure}[t]{3cm}
		\centering	
		\begin{tikzpicture}[scale=2]
		\draw[] (0,0)--(0,.5)(0,1)--(0,1.5) (.60,0)--(.60,.5)(.60,1)--(.60,1.5);
		\draw[dotted](0,.5)--(0,1) (.60,.5)--(.60,1);
		
		\filldraw[black] (0,0) circle (.5pt) node[anchor=north]  {$a_1$}; 
		\filldraw[black] (.6,0) circle (.5pt) node[anchor=north] {$b_1$}; 
		\filldraw[black] (0,1.5) circle (.5pt) node[anchor=south] {$a_m$}; 
		\filldraw[black] (.60,1.5) circle (.5pt) node[anchor=south] {$b_n$}; 
		\end{tikzpicture}
		\caption{}\label{fig15}
	\end{subfigure}
	\quad
	\begin{subfigure}[t]{3cm}
		\centering	
		\begin{tikzpicture}[scale=2]
		\draw[] (0,0)--(0,.5)(0,1)--(0,1.5) (.60,0)--(.60,.5)(.60,1)--(.60,1.5);
		\draw[dotted] (0,0)--(.60,.75) (0,.5)--(0,1) (.60,.5)--(.60,1) ;
		
		\filldraw[black] (0,0) circle (.5pt) node[anchor=north]  {$a_1$}; 
		\filldraw[black] (.6,0) circle (.5pt) node[anchor=north] {$b_1$}; 
		\filldraw[black] (0,1.5) circle (.5pt) node[anchor=south] {$a_m$}; 
		\filldraw[black] (.60,1.5) circle (.5pt) node[anchor=south] {$b_n$};  
		\end{tikzpicture}
		\caption{}\label{fig16}
	\end{subfigure}
	\quad
	\begin{subfigure}[t]{3cm}
		\centering	
		\begin{tikzpicture}[scale=2]
		\draw[] (0,0)--(0,.5)(0,1)--(0,1.5) (.60,0)--(.60,.5)(.60,1)--(.60,1.5);
		\draw[dotted] (0,.75)--(.60,1.5)(0,.5)--(0,1) (.60,.5)--(.60,1);
		\filldraw[black] (0,0) circle (.5pt) node[anchor=north]  {$a_1$}; 
		\filldraw[black] (.6,0) circle (.5pt) node[anchor=north] {$b_1$}; 
		\filldraw[black] (0,1.5) circle (.5pt) node[anchor=south] {$a_m$}; 
		\filldraw[black] (.60,1.5) circle (.5pt) node[anchor=south] {$b_n$}; 
		\end{tikzpicture}
		\caption{}\label{fig17}
	\end{subfigure}
	\quad
	\begin{subfigure}[t]{3cm}
		\centering	
		\begin{tikzpicture}[scale=2]
		\draw[] (0,0)--(0,.5)(0,1)--(0,1.5) (.60,0)--(.60,.5)(.60,1)--(.60,1.5);
		\draw[dotted] (0,.75)--(.60,1.5)  (0,0)--(.60,.75)(0,.5)--(0,1) (.60,.5)--(.60,1);
		
		\filldraw[black] (0,0) circle (.5pt) node[anchor=north]  {$a_1$}; 
		\filldraw[black] (.6,0) circle (.5pt) node[anchor=north] {$b_1$}; 
		\filldraw[black] (0,1.5) circle (.5pt) node[anchor=south] {$a_m$}; 
		\filldraw[black] (.60,1.5) circle (.5pt) node[anchor=south] {$b_n$}; 
		\end{tikzpicture}
		\caption{}\label{fig18}
	\end{subfigure}
	\caption{}\label{fig20}
\end{figure}

\subsection{Graph Theory}\label{subsec:graphtheory}
Let $G$ be a simple graph on the vertex set $[n]$. The {\em clique complex} (or {\em flag complex}) $\Delta(G)$ associated to $G$ is a simplicial complex defined in the following way: $\Delta(G)$ has same vertices as $G$ and the simplices of $\Delta(G)$ are exactly the subsets $F$ of $[n]$ for which every pair in $F$ is an edge of $G$. A graph G is called {\em chordal} if every induced cycle in $G$ of length $\geq 4$ has a chord, i.e.,  there is an edge in $G$ connecting two nonconsecutive vertices of the cycle. Let $\Delta$ be a simplicial complex. The Stanley-Reisner ideal $I_\Delta$ generated by quadratics has linear resolution if and only if $\Delta =\Delta(G)$ for some chordal graph $G$~\cite[Theorem\ 1]{Fro90}.
		
Let $P$ be a poset. The {\em comparability graph} $G_P$ of $P$ is a graph on the underlying set of $P$ such that $\{x,y\}$ is an edge of $G_P$ if and only if $x$ and $y$ are comparable in $P$. Let $\Delta(P)$ be the order complex of $P$. It is known that $\Delta(P) =\Delta(G_P)$.

\subsection{Squarefree divisor complexes}\label{squarefree}
Let  $H\subset \NN^n$ be an affine semigroup and $K[H]$ be the semigroup ring attached to it. Suppose that $h_1,\ldots, h_m\in \NN^n$ is the unique minimal set of generators of $H$. We consider the polynomial ring $T=K[t_1,\ldots,t_n]$ in $n$ variables. Then $K[H]$ is the subring of $T$ generated by the monomials $u_i=\prod_{j=1}^nt_j^{h_i(j)}$ for $1\leq i \leq m$, where $h_i(j)$ denotes the $j$th component of the integer vector $h_i$. Consider a $K$-algebra map $S=K[x_1,\ldots,x_m]\to K[H]$ with $x_i\mapsto u_i$ for $i=1,\ldots,m$.  Let $I_H$ be the kernel of this $K$-algebra map. Set  $\deg x_i=h_i$ to assign a $\ZZ^n$-graded ring structure to $S$. Let $\mm$ be the graded maximal $S$-ideal. Then $K[H]$ as well as $I_H$ become $\ZZ^n$-graded $S$-modules. Thus, $K[H]$ admits a minimal $\ZZ^n$-graded  $S$-resolution $\FF$.
  
Given $h\in H$, we define the {\em squarefree divisor complex} $\Delta_h$ as follows:

 $$\Delta_h:=\{F \subseteq [m] :\prod_{i\in F}^{} u_i \ \text{divides}\ t_1^{h(1)}\cdots t_n^{h(n)}\ \text{in}\ K[H]\}.$$

 Clearly, $\Delta_h$ is a simplicial complex. We denote the $i^{\text{th}}$ reduced simplicial homology of a simplicial complex $\Delta$ with coefficients in $K$ by $\widetilde{H}_{i}(\Delta, K)$.
  
  \begin{Proposition}\cite[Proposition\ 1.1]{BH}, \cite[Theorem\ 12.12]{[STURM96]}
  	\label{bh}
  	With the notation and assumptions introduced one has  $\Tor_i(K[H],K)_h\iso\widetilde{H}_{i-1}(\Delta_h, K)$. In particular,
  	$$\beta_{ih}(K[H])=\dim_K\widetilde{H}_{i-1}(\Delta_h, K).$$
  \end{Proposition}
  
   Let $H\subset \NN^n$ be an affine semigroup generated by $h_1,\ldots, h_m$. An affine subsemigroup $H'\subset H$ generated by  a subset of $\{h_1,\ldots, h_m\}$ is called a {\em homologically pure} subsemigroup of $H$ if for all $h\in H'$ and all $h_i$ with $h-h_i\in H$,  it follows that $h_i\in H'$. Let $H'$ be generated by a subset $\mathcal{X}$ of $\{h_1,\ldots,h_m\}$, and let $S'=K[\{x_i : h_i\in \mathcal{X}\}]\subseteq S$. Therefore, $K[H']$ has $\ZZ^{n}$-graded minimal free $S'$-resolution. Let $\Syz^S_i(K[H])$ denotes the $i^{\text{th}}$ syzygy module of $K[H]$ over $S$.
    
 We need the following proposition several times in this paper. 
  
\begin{Proposition}\cite[Corollary\ 3.4]{[EHH15]}
  	\label{homologicallypure}
 Let $H'$ be a homologically pure subsemigroup  of $H$. Then, any minimal set of generators of $\Syz^{S'}_i(K[H'])$ is part of a minimal set of generators  of  $\Syz^S_i(K[H])$ for all $i$. Moreover, $\beta^{S'}_{ij}(K[H'])\leq \beta^{S}_{ij}(K[H])$ for all $i$ and $j$.
\end{Proposition}
  
Let $L=\MI(P)$ be a distributive lattice  with $P=\{p_1,\ldots,p_n\}$. For $ \alpha \in L$, define a $(n+1)$-tuple $h_\alpha$ such that for $1 \leq i \leq n$, 
 \[
 \begin{cases} 
 1 &  \text{at $1^{st}$ position}  \text{,}  \\
 1 & \text{at $(i+1)^{th}$ position if} \ p_i \in \alpha  \text{,}  \\
 0 & \text{at $(i+1)^{th}$ position if} \ p_i \notin \alpha.
 \end{cases}
 \]
 Let $H$ be the affine semigroup generated by $\{h_\alpha : \alpha \in L\}$. Then, we have $K[H]=R[L]$.
 
\section{Homologically pure subsemigroups}\label{sec:homologicallypure}

In this section, we identify two kinds of homologically pure subsemigroups of an affine semigroup associated to a Hibi ring and we use them to conclude results about property $N_p$ of Hibi rings. The first one is the following and the second one is in Notation~\ref{construction1}.

Let $L=\MI(P)$ be a distributive lattice. Let $\beta, \gamma\in L$ such that $\beta \leq \gamma$. Define $L_1 = \{\alpha \in L : \beta \leq \ \alpha \leq\gamma\}$. Clearly, $L_1$ is a sublattice of $L$. Let $H$ be the affine semigroup associated to $L$ and let $H_1$ be the affine subsemigroup of $H$ generated by $\{h_\alpha : \alpha \in L_1\}$.

\begin{Proposition}\label{semigroup}
	Let $H$ and $H_1$ be as defined above. Then $H_1$ is a  homologically pure subsemigroup of $H$.
\end{Proposition}
\begin{proof}
	We show that if $\alpha \notin L_1$, then $h - h_\alpha \notin H$ for all  $h \in H_1$. Suppose $\alpha \notin L_1$ then either $\alpha \nleq \gamma$ or $\alpha \ngeq \beta$.\par
	
	If $\alpha \nleq \gamma$, then there exists a $p_i \in \alpha$ such that $p_i \notin \gamma$. So $i^{th}$ entry of $h_\alpha$ is 1 but for any $\alpha' \in L_1$, $i^{th}$ entry of $h_{\alpha'}$ is 0. Hence, $h - h_\alpha \notin H$ for all  $h \in H_1$.\par
	
	If $\alpha \ngeq \beta$, then there exists a $p_j \in \beta$ such that $p_j \notin \alpha$. So $(j+1)^{th}$ entry of $h_{\alpha}$ is 0 but for any $\alpha' \in L_1$, $(j+1)^{th}$ entry of $h_{\alpha'}$ is 1. Therefore, for any $h\in H_1$, the first and $(j+1)^{th}$ entries of $h$ have same value, say $r_h$. Thus, for all $h \in H_1$, the first entry and $(j+1)^{th}$ entry of $h - h_\alpha$ are $r_h-1$ and $r_h$ respectively. Therefore, $h - h_\alpha \notin H$ for all  $h \in H_1$.
\end{proof}

\begin{Proposition}\label{sublattice}
	Let $L$ and $L_1$ be as above. Let $\beta = \{p_{a_1},\ldots,p_{a_r}\}$ and $\gamma = \{p_{a_1},\ldots,p_{a_r},\ab p_{b_1},\ldots,p_{b_s}\}$. Then, the induced subposet $P_1$ of $P$ on the set $\{p_{b_1},\ldots,p_{b_s}\}$ is isomorphic to the poset of join-irreducible elements of $L_1$.
\end{Proposition}
\begin{proof}
	The idea of the proof is based on the proof of \cite[Theorem\ 6.4]{[HHO18]}. For any two finite posets $Q$ and $Q'$, if $\MI(Q)\cong \MI(Q')$ then $Q\cong Q'.$ So it is enough to prove that $\MI(P_1) \cong L_1.$  Define a map 
	$$\varphi: \MI(P_1) \to L_1$$ by \[
	\varphi(\alpha) = (\mathop{\vee}_{i=1}^{r} p_{a_i}) \vee (\mathop{\vee}_{p\in \alpha} p).
	\]
	In particular, $\varphi(\emptyset) = \mathop{\vee}_{i=1}^{r} p_{a_i}$. Clearly, $\varphi$ is order-preserving.\par
	Let $\alpha$ and $\delta$ be two order ideals of $P_1$ with $\alpha \neq \delta$, say $\delta \nleq \alpha$. Let $p_0$ be a maximal element of $\delta$ with $p_0 \notin \alpha$. We show that $\varphi(\alpha) \neq \varphi(\delta)$. Suppose, on the contrary, $\varphi(\alpha) = \varphi(\delta)$, then $$(\mathop{\vee}_{i=1}^{r} p_{a_i}) \vee (\mathop{\vee}_{p\in \alpha} p) = (\mathop{\vee}_{i=1}^{r} p_{a_i}) \vee (\mathop{\vee}_{q\in \beta} q).$$ 
	Since $L_1$ is distributive, it follows that
	$$((\mathop{\vee}_{i=1}^{r} p_{a_i}) \vee (\mathop{\vee}_{p\in \alpha} p))\wedge p_0 = (\mathop{\vee}_{i=1}^{r} (p_{a_i}\wedge p_0)) \vee (\mathop{\vee}_{p\in \alpha} (p\wedge p_0)).$$
	Since $p_0$ is join-irreducible and for any $p\in P$, $p\wedge p_0\leq p_0$. It follows that $(\mathop{\vee}_{i=1}^{r} p_{a_i} \vee (\mathop{\vee}_{p\in \alpha} p))\wedge p_0 < p_0$. However, since $p_0 \in \delta$,  $(\mathop{\vee}_{i=1}^{r} p_{a_i} \vee (\mathop{\vee}_{q\in \beta} q)) \wedge p_0 = p_0$. This is a contradiction. Hence, $\varphi$ is injective.\par 
	Since each $a\in L_1$ can be the join of the join-irreducible elements $p$ with $p \leq a$ in $L_1$, it follows that $\varphi(\alpha) =a$, where $\alpha$ is an order ideal of $P_1$ consisting of those $p \in P_1$ with $p\leq a$. Thus, $\varphi$ is surjective.\par
	
	Now, $\varphi^{-1}: L_1 \to \MI(P_1)$ is  defined as follows: for $x\in L_1$, $$\varphi^{-1}(x)= \{p\in L_1: p\leq x, p \ \text{is a join-irreducible}\}\setminus \mathop{\cup}_{i=1}^{r} p_{a_i}.$$
	Clearly, $\varphi^{-1}$ is order-preserving. Hence the proof.
\end{proof}

We now try to understand how we are going to use the above propositions. Suppose for a distributive lattice $\mathcal{L}$,  we want to prove that  $\beta_{ij}(R[\mathcal{L}])\neq 0$ for some $i,j$. The idea of the proof is to reduce the lattice $\mathcal{L}$ to a suitably chosen sublattice $\mathcal{L}_1$. Proposition~\ref{sublattice}  describes the subposet of join-irreducibles of $\mathcal{L}_1$. Then, by Propositions~\ref{semigroup} and \ref{homologicallypure}, if $\beta_{ij}(R[\mathcal{L}_1])\neq 0$, then $\beta_{ij}(R[\mathcal{L}])\neq 0$. More precisely,

\begin{discussionbox}\label{subsemigroup}\normalfont
	Let $P$ be a poset. Let $B$ and $ B'$ be two antichains of $P$ such that for each $p\in B$ there is a $q\in B'$ such that $p\leq q$ and for each $q'\in B'$ there is a $p'\in B$ such that $p'\leq q'$. Furthermore, let $\gamma= \{p\in P : p\leq q\ \text{for some}\ q\in B'\}$ and $\beta' = \{p\in P:p' \leq p\leq q\ \text{for some}\ p'\in B, q\in B' \}$. Let $\beta = \gamma \setminus \beta'$. Note that $\beta,\gamma$ are the order ideals of $\MI(P)$ with $\beta < \gamma$. Let $L_1 = \{\alpha \in \MI(P) : \beta \leq \ \alpha \leq\gamma\}$. Furthermore, let $H_1$ be the affine subsemigroup of $H$ generated by $\{h_\alpha : \alpha \in L_1\}$. Then, by Proposition~\ref{semigroup}, $H_1$ is a  homologically pure subsemigroup of $H$. Also, by Proposition~\ref{sublattice}, the induced subposet $P_1$  of $P$ on the set $\gamma \setminus \beta$ is isomorphic to the poset of join-irreducible elements of $L_1$.
	Furthermore, by Proposition~\ref{homologicallypure}, $\beta_{ij}(R[L_1]) \leq \beta_{ij}(R[L])$ for all $i, j$.
\end{discussionbox}

\begin{Example}\normalfont
    In this example, we illustrate the construction in the above discussion. 
    Let $P$ be as shown in Figure~\ref{fig:poset}. 
    Then, $\MI(P)$ is as shown in Figure~\ref{fig:lattice}.
    Under the notations of Discussion~\ref{subsemigroup}, let $B=\{p_3,p_4\}$ and $B'=\{p_6,p_7\}$.
    Then, $\gamma = P\setminus \{p_8\}$ and $\beta' = \{p_3,p_4,p_5,p_6,p_7\}$.
    Thus, $\beta = \gamma\setminus\beta' = \{p_1,p_2\}$. 
\end{Example}

	\begin{figure}[h]
		\begin{subfigure}[t]{8cm}
			\centering	
			\begin{tikzpicture}[scale=1]
                \draw[]   (0,0)--(0,3) (1,0)--(1,3)
                          (0,3)--(1,1) (0,1)--(1,2)  
	                    	;	
			\filldraw[black] (0,0) circle (.5pt) node[anchor=north]  {$p_1$}; 
			\filldraw[black] (0,1) circle (.5pt) node[anchor=east] {$p_3$}; 
			\filldraw[black] (1,0) circle (.5pt) node[anchor=north] {$p_2$};  
			\filldraw[black] (1,1) circle (.5pt) node[anchor=west]  {$p_4$}; 
			\filldraw[black] (0,2) circle (.5pt) node[anchor=east] {$p_5$}; 
			\filldraw[black] (1,2) circle (.5pt) node[anchor=west] {$p_6$}; 
			\filldraw[black] (0,3) circle (.5pt) node[anchor=south] {$p_7$}; 
			\filldraw[black](1,3)  circle (.5pt) node[anchor=south]{$p_8$};  	 		 		
			\end{tikzpicture}
			\caption{}\label{fig:poset}
		\end{subfigure}
		\quad
		\begin{subfigure}[t]{7cm}
			\centering	
			\begin{tikzpicture}[scale=1]
                \draw[]      (0,0)--(0,3) (1,0)--(1,3) (2,0)--(2,4) (3,2)--(3,4) (4,2)--(4,4)
                             (0,0)--(2,0) (0,1)--(2,1) (0,2)--(4,2) (0,3)--(4,3) (2,4)--(4,4);

			\filldraw[black] (0,0) circle (.5pt) node[anchor=north]  {$\emptyset$}; 
            \filldraw[black] (0,1) circle (.5pt) node[anchor=east] {$\{p_1\}$}; 
            \filldraw[black] (1,0) circle (.5pt) node[anchor=north] {$\{p_2\}$}; 
            \filldraw[black] (0,2) circle (.5pt) node[anchor=east]  {$\{p_1,p_3\}$}; 
            \filldraw[black] (0,3) circle (.5pt) node[anchor=east] {$\{p_1,p_3,p_5\}$}; 
            \filldraw[black] (2,0) circle (.5pt) node[anchor=west] {$\{p_2,p_4\}$}; 
            \filldraw[black] (2,4) circle (.5pt) node[anchor=east] {$\{p_1,p_2,p_3,p_4,p_5,p_7\}$}; 
            
            \filldraw[black] (1,1) circle (1pt) node[anchor= south east] {$\beta$}; 
            \filldraw[black] (3,4) circle (1pt) node[anchor=south] {$\gamma$}; 
			\end{tikzpicture}
			\caption{}\label{fig:lattice}
		\end{subfigure}
		\caption{}
	\end{figure} 

\begin{notationbox}\label{construction1}\normalfont
	For a poset $P$, let $X_P$ and $Y_P$ be the sets of minimal and maximal elements of $P$ respectively. Define $X'_P = \{q\in P : p\lessdot q \ \text{for some} \ p\in X_P\}$ and $Y'_P = \{p\in P : p \lessdot q \ \text{for some} \ q\in Y_P\}$. When the context is clear, we will omit the subscripts and denote $X_P, X'_P,Y_P$ and $Y'_P$ by $X, X',Y$ and $Y'$ respectively.\par
	
	Let $P$ be a poset. For $x,y \in P$ with $x<y$, define $L_1 := \{\alpha \in\MI(P) : \text{if} \ x \in \alpha\  \text{then} \ y \in \alpha \}$. 
    It is easy to see that $L_1$ is a sublattice of $\MI(P)$. 
    Let $P_1$ be the poset on the set $P \setminus \{p\in P: x\leq p <y\}$ defined by the following order relations: if $p,q \in P_1,$ then $p\leq q$ in $P_1$ if  either\par 
	(1) $p\in P_1\setminus{\{y\}}, q \in P_1$ and $p \leq q$ in $P$ or \par 
	(2) $p =y$ and there is a $p' \in \{a\in P: x\leq a \leq y\}$ such that $ p' \leq q$ in $P$.\par
	Let $H$ be the semigroup corresponding to $\MI(P)$ and $H'$ be the subsemigroup of $H$ corresponding to $L_1$.
\end{notationbox}

\begin{Lemma}\label{reduceposet}
	Let $P, P_1, L_1, H, H'$ be as in Notation~\ref{construction1}. Then $L_1 \cong \MI(P_1)$ and $H'$ is a homologically pure subsemigroup of $H$.
\end{Lemma}

\begin{proof}
	Define a map 
	$$\varphi: \MI(P_1) \to L_1$$ by
	\[ \varphi(\alpha) =  \begin{cases} 
	\alpha & \text{if}\ y\notin\alpha, \\
	\alpha \cup \{p\in P: x\leq p <y\} & \text{if} \ y \in \alpha.  \\
	\end{cases}
	\]
	Clearly, $\varphi$ is order-preserving. If $\gamma \in L_1$, then $\varphi(\gamma')=\gamma$, where $\gamma'= \gamma\setminus\{p\in P: x\leq p <y\}$. 
    Hence, $\varphi$ is surjective. Now we claim that for any $\alpha\in \MI(P_1)$, $\varphi(\alpha) \cap P_1 =\alpha$. If $y\in \alpha$, 
    then $\varphi(\alpha) \cap P_1= (\alpha \cup \{p\in P: x\leq p <y\}) \cap P_1= \alpha$ and  if $y\notin \alpha$, then $\varphi(\alpha)=\alpha$. Therefore,  if $\varphi(\alpha)=\varphi(\beta)$ for any $\alpha,\beta \in \MI(P_1)$ then $\alpha= \beta$. This proves that $\varphi$ is injective. \par
	
	Now, $\varphi^{-1}: L_1 \to \MI(P_1)$ is  defined as follows: for $a\in L_1$, $$\varphi^{-1}(a)= \{p\in L_1: p\leq a, p \ \text{is a join-irreducible}\}\setminus \{p\in P: x\leq p <y\}.$$
	Clearly, $\varphi^{-1}$ is order-preserving. Hence, $\varphi$ is an isomorphism.	
	
	To  prove that $H'$ is a homologically pure subsemigroup of $H$, we show that if $\alpha \notin L_1$ then $h - h_\alpha \notin H$ for all  $h \in H'$. Suppose $\alpha \notin L_1$ then  $x\in\alpha$ but $y\notin\alpha$. Let $ h = \sum_{i =1}^{s} h_{\beta_i} \in H'$ and let the position corresponding to $x$ of $h$ be $r$. Then the positions corresponding to $x$ and $y$ of $h-h_{\alpha}$ are $r-1$ and $r$ respectively. Hence, $h-h_{\alpha} \notin H$.
\end{proof}

\begin{discussionbox}\label{reducingposet}\normalfont
	For a poset $P_0$, let $X_{P_{0}},$ $Y_{P_{0}},$ $X'_{P_{0}}$ and $Y'_{P_{0}}$ be as defined in Notation~\ref{construction1}. If there is an $x\in X'_{P_{0}}$ and a $y\in Y'_{P_{0}}$ with $x<y$, reduce $P_0$ to $P_1$, using the methods in Notation~\ref{construction1}. Observe that $y\in X'_{P_{1}}\cap Y'_{P_{1}}$, $X_{P_{0}}=X_{P_{1}}, Y_{P_{0}}=Y_{P_{1}}$ and $\#P_{1} =  \#(P_0\setminus\{p\in P_0: x\leq p <y\}) \leq \#P_0-1$. Repeating it, we get a sequence of posets $P_0,\ldots, P_n$, where $n\leq \#P_0-\#X_0-\#Y_0-1$ such that for each $0\leq i\leq n-1,$ there is an $x\in X'_{P_{i}}$ and $y\in Y'_{P_{i}}$ with $x<y$ and $P_{i}$ is reduced to $P_{i+1}$ as in Notation~\ref{construction1}. Moreover, there is no $x\in X'_{P_{n}}$ and $y\in Y'_{P_{n}}$ with the property $x<y$. Here, $P_n$ is a  poset defined on the set $X_{P_{0}} \cup Y_{P_{0}} \cup Y'_{P_{n}}$ and $\rank(P_n)\leq2$. An example of this reduction is given in Figure~\ref{fig8}. By Lemma~\ref{reduceposet} and Proposition~\ref{homologicallypure}, if $\beta_{24}(R[\MI(P_i)]) \neq 0$ for some $1\leq i \leq n$, then $\beta_{24}(R[\MI(P_0)]) \neq 0$.
\end{discussionbox}

	\begin{figure}[h!]
		\begin{subfigure}[t]{7cm}
			\centering	
			\begin{tikzpicture}[scale=2]
			\draw[] (.8,0)--(.4,.6)--(0,0)--(-.4,.6)--(-.4,1.1)--(-.8,1.5) (-.4,1.1)--(0,1.5) (.4,.6)--(.4,1.1)
			(.8,0)--(1,.6)--(1.3,1.1)--(1.3,1.5)  (1.5,0)--(1.5,.6)--(1.3,1.1);
			
			\filldraw[black] (0,0) circle (.5pt) node[anchor=north]  {$p_1$}; 
			\filldraw[black] (1.5,0) circle (.5pt) node[anchor=north] {$p_3$}; 
			\filldraw[black] (.8,0) circle (.5pt) node[anchor=north] {$p_2$}; 
			\filldraw[black] (-.4,.6) circle (.5pt) node[anchor=east] {$p_4$}; 
			\filldraw[black] (.4,.6) circle (.5pt) node[anchor=east]  {$p_5$}; 
			\filldraw[black] (1,.6) circle (.5pt) node[anchor=east] {$p_6$}; 
			\filldraw[black] (1.5,.6) circle (.5pt) node[anchor=west] {$p_7$}; 
			\filldraw[black] (-.4,1.1) circle (.5pt) node[anchor=east] {$p_8$}; 
			\filldraw[black] (.4,1.1) circle (.5pt) node[anchor=east]  {$p_9$}; 
			\filldraw[black] (1.3,1.1) circle (.5pt) node[anchor=west] {$p_{10}$}; 
			\filldraw[black] (0,1.5) circle (.5pt) node[anchor=south] {$p_{11}$}; 
			\filldraw[black] (-.8,1.5) circle (.5pt) node[anchor=south] {$p_{12}$};  
			\filldraw[black] (1.3,1.5) circle (.5pt) node[anchor=south] {$p_{13}$}; 	
			\end{tikzpicture}
			\caption{}\ 
		\end{subfigure}
		\quad
		\begin{subfigure}[t]{6cm}
			\centering	
			\begin{tikzpicture}[scale=2]
			\draw[] (.8,0)--(.4,.6)--(0,0)--(-.4,.6)--(-.8,1.1) (-.4,.6)--(0,1.1) (.4,.6)--(.4,1.1) 
			(.8,0)--(1,.6)--(1.3,1.1)--(1.3,1.5)  (1.5,0)--(1.5,.6)--(1.3,1.1);
			\filldraw[black] (0,0) circle (.5pt) node[anchor=north]  {$p_1$}; 
			\filldraw[black] (1.5,0) circle (.5pt) node[anchor=north] {$p_3$}; 
			\filldraw[black] (.8,0) circle (.5pt) node[anchor=north] {$p_2$}; 
			\filldraw[black] (-.4,.6) circle (.5pt) node[anchor=east] {$p_8$}; 
			\filldraw[black] (.4,.6) circle (.5pt) node[anchor=east]  {$p_5$}; 
			\filldraw[black] (1,.6) circle (.5pt) node[anchor=east] {$p_6$}; 
			\filldraw[black] (1.5,.6) circle (.5pt) node[anchor=west] {$p_7$};  
			\filldraw[black] (.4,1.1) circle (.5pt) node[anchor=south]  {$p_{9}$}; 
			\filldraw[black] (1.3,1.1) circle (.5pt) node[anchor=west] {$p_{10}$}; 
			\filldraw[black] (1.3,1.5) circle (.5pt) node[anchor=south] {$p_{13}$};  
			\filldraw[black] (-.8,1.1) circle (.5pt) node[anchor=south] {$p_{11}$}; 
			\filldraw[black] (0,1.1) circle (.5pt) node[anchor=south] {$p_{12}$}; 
			\end{tikzpicture}
			\caption{If $p_4 \in \alpha$, then $p_8\in \alpha$}\ 
		\end{subfigure} 
		\quad
		\begin{subfigure}[t]{7cm}
			\centering	
			\begin{tikzpicture}[scale=2]
			\draw[] (.8,0)--(.4,.6)--(0,0)--(-.4,.6)--(-.8,1.1) (-.4,.6)--(0,1.1) (.4,.6)--(.4,1.1) 
			(.8,0)--(1,.6)--(1.3,1.1)  (1.5,0)--(1.5,.3)--(1,.6);
			\filldraw[black] (0,0) circle (.5pt) node[anchor=north]  {$p_1$}; 
			\filldraw[black] (1.5,0) circle (.5pt) node[anchor=north] {$p_3$}; 
			\filldraw[black] (.8,0) circle (.5pt) node[anchor=north] {$p_2$}; 
			\filldraw[black] (-.4,.6) circle (.5pt) node[anchor=east] {$p_8$}; 
			\filldraw[black] (.4,.6) circle (.5pt) node[anchor=east]  {$p_5$}; 
			\filldraw[black] (1,.6) circle (.5pt) node[anchor=east] {$p_{10}$};  
			\filldraw[black] (.4,1.1) circle (.5pt) node[anchor=south]  {$p_{9}$}; 
			\filldraw[black] (1.3,1.1) circle (.5pt) node[anchor=south] {$p_{13}$};   
			\filldraw[black] (-.8,1.1) circle (.5pt) node[anchor=south] {$p_{11}$}; 
			\filldraw[black] (0,1.1) circle (.5pt) node[anchor=south] {$p_{12}$}; 
			\filldraw[black] (1.5,0) circle (.5pt) node[anchor=north] {$p_3$}; 
			\filldraw[black] (1.5,.3) circle (.5pt) node[anchor=west] {$p_7$}; 
			
			\end{tikzpicture}
			\caption{If $p_6 \in \alpha$, then $p_{10}\in \alpha$}\ 
		\end{subfigure} 
		\quad
		\begin{subfigure}[t]{8cm}
			\centering	
			\begin{tikzpicture}[scale=2]
			\draw[] (.8,0)--(.4,.6)--(0,0)--(-.4,.6)--(-.8,1.1) (-.4,.6)--(0,1.1) (.4,.6)--(.4,1.1) 
			(.8,0)--(1,.6)--(1,1.1)  (1.4,0)--(1,.6);
			\filldraw[black] (0,0) circle (.5pt) node[anchor=north]  {$p_1$}; 
			\filldraw[black] (1.4,0) circle (.5pt) node[anchor=north] {$p_3$}; 
			\filldraw[black] (.8,0) circle (.5pt) node[anchor=north] {$p_2$}; 
			\filldraw[black] (-.4,.6) circle (.5pt) node[anchor=east] {$p_8$}; 
			\filldraw[black] (.4,.6) circle (.5pt) node[anchor=east]  {$p_5$}; 
			\filldraw[black] (1,.6) circle (.5pt) node[anchor=east] {$p_{10}$};  
			\filldraw[black] (.4,1.1) circle (.5pt) node[anchor=south]  {$p_{9}$}; 
			\filldraw[black] (1,1.1) circle (.5pt) node[anchor=south] {$p_{13}$};   
			\filldraw[black] (-.8,1.1) circle (.5pt) node[anchor=south] {$p_{11}$}; 
			\filldraw[black] (0,1.1) circle (.5pt) node[anchor=south] {$p_{12}$}; 
			\end{tikzpicture}
			\caption{If $p_7 \in \alpha$, then $p_{10}\in \alpha$}\ 
		\end{subfigure}  	      	 	   	
		\caption{}\label{fig8}
	\end{figure}

\begin{Example}\normalfont
	In this example, we show that the converse of the conclusion in Discussion~\ref{reducingposet} may not be true. Let $P$ be a poset as shown in Figure~\ref{fig19}. By Lemma~\ref{cardinality4betti}, $\beta_{24}(R[\MI(P)]) \neq 0$. Now, let $x= p_4$ and $y=p_7$. 
    Reduce $P$ to $P_1$, using the methods of Notation~\ref{construction1}. 
    Since $p_4 < p_6$ in $P$, we have $p_7< p_6$ in $P_1$ by the definition of $P_1$, and the order relations of $P$ are also order relations of $P_1$.
    Thus, $P_1$ is as shown in Figure~\ref{fig21}. By Theorem~\ref{ene}, $\beta_{24}(R[\MI(P_1)]) = 0$. 
	
	\begin{figure}[h]
		\begin{subfigure}[t]{8cm}
			\centering	
			\begin{tikzpicture}[scale=2]
			\draw[](-.4,0)--(-.4,.8)--(1,0)--(.4,.6)--(.4,1.1)--(1,1.7)--(1.6,1.1)--(1.6,.6)-- (1,0)  (.4,.6)--(1.6,1.1)(.4,1.1)--(1.6,.6);
			
			\filldraw[black] (-0.4,0) circle (.5pt) node[anchor=north]  {$p_1$}; 
			\filldraw[black] (-.4,.8) circle (.5pt) node[anchor=south] {$p_3$}; 
			\filldraw[black] (1,0) circle (.5pt) node[anchor=north] {$p_2$};  
			\filldraw[black] (0.4,.6) circle (.5pt) node[anchor=north]  {$p_4$}; 
			\filldraw[black] (.4,1.1) circle (.5pt) node[anchor=east] {$p_6$}; 
			\filldraw[black] (1,1.7) circle (.5pt) node[anchor=south] {$p_8$}; 
			\filldraw[black] (1.6,1.1) circle (.5pt) node[anchor=west] {$p_7$}; 
			\filldraw[black](1.6,.6)circle(.5pt)node[anchor=north]{$p_5$};  	 		 		
			\end{tikzpicture}
			\caption{}\label{fig19}
		\end{subfigure}
		\quad
		\begin{subfigure}[t]{7cm}
			\centering	
			\begin{tikzpicture}[scale=2]
			\draw[](0,0)--(0,.7)--(.7,0)--(.7,.6)--(.7,1.1)--(.7,1.6)--(.7,2);
			\filldraw[black] (0,0) circle (.5pt) node[anchor=north]  {$p_1$}; 
			\filldraw[black] (0,.7) circle (.5pt) node[anchor=south] {$p_3$}; 
			\filldraw[black] (.7,0) circle (.5pt) node[anchor=north] {$p_2$}; 
			\filldraw[black] (.7,.6) circle (.5pt) node[anchor=west]  {$p_5$}; 
			\filldraw[black] (.7,1.1) circle (.5pt) node[anchor=west] {$p_7$}; 
			\filldraw[black] (.7,1.6) circle (.5pt) node[anchor=west] {$p_6$}; 
			\filldraw[black] (.7,2) circle (.5pt) node[anchor=west] {$p_8$}; 
			\end{tikzpicture}
			\caption{}\label{fig21}
		\end{subfigure}
		\caption{}
	\end{figure} 
\end{Example}

\section{Property $N_2$ of Hibi rings}\label{propertyn2}

In this section, we prove some sufficient conditions regarding when Hibi rings do not satisfy  property $N_2$. 

\subsection{}

Here, we prove Theorem~\ref{planar}. It shows how to reduce checking property $N_2$ to a planar distributive sublattice. We begin by proving some  relevant lemmas.

\begin{Lemma}\label{bettitensor}\cite[Problem\ 2.16]{[HHO18]}
Let $K$ be a field, $S= K[x_1 ,\ldots, x_n]$ and $T = K[y_1 ,\ldots,y_m]$ be two polynomial rings. Let $M$ be a finitely generated graded $S$-module
and $N$ be a finitely generated graded $T$-module. Then $M\tensor_{K}N$ is a finitely generated graded $S\tensor_{K}T$-module and
$$\beta_{pq}(M\tensor_{K}N) = \sum\beta_{p_{1}q_{1}}(M)\beta_{p_{2}q_{2}}(N),$$
where the sum is taken over all $p_1$ and $p_2$ with $p_1 + p_2 = p$, and over all $q_1$ and $q_2$ with $q_1 + q_2 = q$.
\end{Lemma}
 
\begin{Lemma}\label{nonsimple}
Let $P$  be a poset and	$p$ be an element of $P$ which is comparable to every element of $P$. Let $P_1=\{q\in P: q<p\}$ and $P_2=\{q\in P: q>p\}$ be induced subposets of $P$. If  $P_1$ and $P_2$ are not chains, then $R[\MI(P)]$ does not satisfy property $N_2$. 
\end{Lemma}

\begin{proof}
 Since $P_1$ and $P_2$ are not chains, $R[\MI(P_1)]$ and $R[\MI(P_2)]$ are not polynomial rings. Therefore, $\beta_{12}(R[\MI(P_{i})])\neq0$ for $i=1,2$. Note that $P$ is the ordinal sum  $P_1\dirsum\{p\}\dirsum P_2$. By Lemma~\ref{ordinalsum}, $R[\MI(P)]=R[\MI(P_1)]\tensor R[\MI(P_2)]$. Hence, $\beta_{24}(R[\MI(P)])\neq0$ by Lemma~\ref{bettitensor}.
\end{proof} 

In \cite{[ENE15]}, Ene proved the above lemma for the case when $\MI(P)$ is a planar distributive lattice.

\begin{Lemma}\label{addingrelations}
Let $P$ be a simple poset such that $\#P = m+n$. Let $\MI(P)$ be a planar distributive lattice such that $\MI(P) \subseteq [(0,0),(m,n)]$ with $m,n\geq 2$. On the underlying set of $P$, let $P'$ be a poset such that every order relation in $P$ is also an order relation in $P'$. Assume that the set of minimal (respectively maximal) elements of $P'$ coincide with the set of minimal (respectively maximal) elements of $P$. If $\beta_{24}(R[\MI(P)]) \neq 0$, then $\beta_{24}(R[\MI(P')]) \neq 0$.
\end{Lemma}

\begin{proof}
	If $P'$ is not simple, then there exists an element $p\in P'$ which is comparable to every element of $P'$. 
    Since $P$ is simple and $\MI(P)$ is a planar distributive lattice, $P$ has exactly two minimal elements and exactly two maximal elements. 
    Since the set of minimal (respectively maximal) elements of $P'$ coincide with the set of minimal (respectively maximal) elements of $P$, we get that $p$ is neither a minimal element nor a maximal element in $P'$. 
    Let $P_1=\{q\in P: q<p\}$ and $P_2=\{q\in P: q>p\}$. Since $P_1$ and $P_2$ are not chains, $\beta_{24}(R[\MI(P')]) \neq 0$ by Lemma~\ref{nonsimple}. So we may assume that $P'$ is simple. Suppose, on the contrary, $\beta_{24}(R[\MI(P')]) = 0$. So the conditions $(i)$ and $(ii)$ of Theorem~\ref{ene} hold for $\MI(P')$. Since $\MI(P')\subseteq \MI(P)$, the conditions $(i)$ and $(ii)$ of Theorem~\ref{ene} also hold for $\MI(P)$ which is a contradiction. Hence the proof.
\end{proof}

 \begin{figure}[h]
 	\begin{subfigure}[t]{7cm}
 		\centering	
 		\begin{tikzpicture}[scale=2]
 		\draw[] (0,0)--(0,1)--(1,0)--(1,1)--(0,0);
 		
 		\filldraw[black] (0,0) circle (.5pt) node[anchor=north]  {$p_1$}; 
 		\filldraw[black] (0,1) circle (.5pt) node[anchor=south] {$p_3$}; 
 		\filldraw[black] (1,0) circle (.5pt) node[anchor=north] {$p_2$}; 
 		\filldraw[black] (1,1) circle (.5pt) node[anchor=south] {$p_4$}; 
 		\end{tikzpicture}
 		\caption{}\label{fig4}
 	\end{subfigure}
 	\quad
 	\begin{subfigure}[t]{7cm}
 		\centering	
 		\begin{tikzpicture}[scale=2]
 		\draw[] (.5,.35)--(1,0) (0,0)--(.5,.35) (.5,.8)--(0,1.25) (.5,.8)--(1,1.25);
 		\draw[dotted](.5,.4)--(.5,.8);
 		\filldraw[black] (0,0) circle (.5pt) node[anchor=north]  {$p_1$}; 
 		\filldraw[black] (0,1.25) circle (.5pt) node[anchor=south] {$p_3$}; 
 		\filldraw[black] (1,0) circle (.5pt) node[anchor=north] {$p_2$}; 
 		\filldraw[black] (1,1.25) circle (.5pt) node[anchor=south] {$p_4$}; 
 		\filldraw[black] (.5,0.35) circle (.5pt) node[anchor=north]  {$q_1$}; 
 		\filldraw[black] (.5,.8) circle (.5pt) node[anchor=west]  {$q_n$}; 
 		\end{tikzpicture}
 		\caption{$n\geq 1$}\label{fig5}
 	\end{subfigure}
 	\caption{}
 \end{figure}

\begin{Lemma}\label{cardinality4betti}
	Let $P$ be a poset such that the poset $P' = \{p_1,...,p_4\}$ of Figure~\ref{fig4}
	is a cover-preserving subposet of $P$. Then $R[\MI(P)]$ does not satisfy property $N_2$.
\end{Lemma}

\begin{proof} 
	Observe that by Theorem~\ref{ene}, $\beta_{24}(R[\MI(P')]) \neq 0$. Let $B=\{p_1,p_2\}$, $B'=\{p_3,p_4\}$. By Discussion~\ref{subsemigroup}, we may replace $P$ by $P_1$,  where $P_1$ is as defined in Discussion~\ref{subsemigroup}, and assume that the sets of minimal and maximal elements of $P$ coincide with the sets of minimal and maximal elements of $P'$ respectively.\\
    Now, suppose that there exists an element $p\in P$ such that $p\notin P'$. Then, we have $p_i<p<p_j$ for some $i\in \{1,2\}$ and $j \in \{3,4\}$. This contradicts that $p_i\lessdot p_j$. Therefore, $P =P'$. This completes the proof.
\end{proof}

\begin{discussionbox}\label{reducesubposet}	
\normalfont Let $P$ be a poset. For $k\geq 1$, let $\mathcal{S}=\mathop{\cup}_{i=1}^{k} \{p_{i,1},\ldots ,p_{i,n_i}\}$ be a subset of the underlying set of $P$. Assume that $\{p_{1,1},\ldots ,p_{k,1}\}$ and $\{p_{1,n_{1}},\ldots ,p_{k,n_k}\}$ are antichains in $P$. Also, assume that for all $1\leq i \leq k$, $\{p_{i,1},\ldots ,p_{i,n_i} \}$  is a chain in $P$ with $p_{i,1}\lessdot\cdots\lessdot p_{i,n_i}$. For $q\in P\setminus \mathcal{S}$, define $\mathcal{S}^P_q := \{p\in \mathcal{S} : q\lessdot p\}$. Let $B=\{p_{1,1},\ldots ,p_{k,1}\}$ and  $B' = \{p_{1,n_{1}},\ldots ,p_{k,n_k}\}$.\par 
Using Discussion~\ref{subsemigroup}, reduce $P$ to $P_1$, where $P_1$ is as defined in Discussion~\ref{subsemigroup}. Let  $x,y\in P_1\setminus \mathcal{S}$ with $x\lessdot y$. Reduce $P_1$ to $P_2$, using the methods of Notation~\ref{construction1}. Observe that $\#P_2=\#P_1-1$, $\mathcal{S}\subset P_2$. Also, $B$ and $B'$ are the sets of minimal and maximal elements of $P_{2}$ respectively. Repeating it, we get a sequence $P_1,\ldots,P_m$, where $m\leq \#P-\#\mathcal{S}$ of posets such that for each $1\leq i\leq m-1$, there exist $x,y\in P_i\setminus \mathcal{S}$ with $x\lessdot y$ and $P_i$ is reduced to $P_{i+1}$ as in Notation~\ref{construction1}. Moreover, there are no $x,y\in P_m\setminus \mathcal{S}$ with the property $x\lessdot y$.\par

Now, we will do more reductions on $P_m$. Let $q\in P_m\setminus \mathcal{S}$ be such that $\#\mathcal{S}^{P_m}_q = 1$, say $\mathcal{S}^{P_m}_q = \{p\}$. We have $q\lessdot p$ in $P_m$. Reduce $P_m$ to $P_{m+1}$, using the methods of Notation~\ref{construction1}. Under this reduction, $\mathcal{S}\subset P_2$, $B$ and $B'$ are the sets of minimal and maximal elements of $P_{m+1}$ respectively. Repeating it,  we get a sequence $P_m,P_{m+1},\ldots,P_s$ of posets such that for each $m\leq i\leq s-1$, there exists a $q\in P_i\setminus \mathcal{S}$ with $\#\mathcal{S}^{P_i}_q = 1$ and $P_i$ is reduced to $P_{i+1}$ as in Notation~\ref{construction1} and there is no $q\in P_s\setminus \mathcal{S}$ with $\#\mathcal{S}^{P_s}_q = 1$.  If $\beta_{ij}(R[\MI(P_l)]) \neq 0$ for some  $i,j$ and $l\in \{1,\ldots, s\}$, then by Discussion~\ref{subsemigroup}, Lemma~\ref{reduceposet} and Proposition~\ref{homologicallypure}, $\beta_{ij}(R[\MI(P)]) \neq 0$.
\end{discussionbox}

\begin{Lemma}\label{rank3simple}
 Let $P$ be a poset and let the poset $P' = \{p_1,...,p_4,q_1,...,q_n\}$ as shown in Figure~\ref{fig5} be a cover-preserving subposet of $P$ for some $n\geq 1$. Then  $R[\MI(P)]$ does not satisfy property $N_2$. 
\end{Lemma}
\begin{proof}
  Note that by Lemma~\ref{nonsimple}, $\beta_{24}(R[\MI(P')]) \neq 0$. Let \[\mathcal{S} = \{p_1,q_1,\ldots,q_n, p_3\} \cup \{p_2,q_1,\ldots,q_n, p_4\}.\] 
  By Discussion~\ref{reducesubposet}, it suffices to show that $R[\MI(P_m)]$ does not satisfy property $N_2$, where $P_m$ is as defined in Discussion~\ref{reducesubposet}. Note that $\{p_1,p_2\}$ and $\{p_3,p_4\}$ are the sets of minimal and maximal elements of $P_m$ respectively. If there exists a cover-preserving subposet of $P_m$ as shown in Figure~\ref{fig4} then $\beta_{24}(R[\MI(P_m)]) \neq 0$. So we may assume that $P_m$ does not contain any cover-preserving subposet as shown in Figure~\ref{fig4}. Let $\mathcal{S}_q$ be as defined in Discussion~\ref{reducesubposet}. There is no $q\in P_m\setminus \mathcal{S}$ with $\mathcal{S}_q = \{p_3,p_4\}$ otherwise $P_m$ will contain a cover-preserving subposet as shown in Figure~\ref{fig4}. So we deduce that $\# \mathcal{S}_q = 1$ for all $q \in P_m\setminus \mathcal{S}$. Now, reduce $P_m$ to $P_s$ as in Discussion~\ref{reducesubposet}. Then $P_s = P'$. This completes the proof.
 \end{proof}

 \begin{Lemma}\label{reflection}
 	Let $(P,\leq)$ be a poset. Then $\MI(P) \iso \MI(P^\partial)$, where $P^\partial$ is the {\em dual poset} of $P$, that is, $(P^\partial,\preceq )$ is the poset with the same underlying set but its order relations are opposite of $P$ i.e. $p\leq q$ if and only if $q\preceq p$. Hence,  $R[\MI(P)] \iso R[\MI(P^\partial)]$.
 \end{Lemma}
 
 \begin{Theorem}\label{planar}
 Let $P$ be a poset. Let $\mathcal{S}=\mathop{\cup}_{i=1}^{2}\{p_{i,1},\ldots ,p_{i,n_i}\}$ be a subset of the underlying set of $P$ such that
 \begin{enumerate}
 	\item
 	for all $1\leq i \leq 2$, $\{p_{i,1},\ldots ,p_{i,n_i} \}$  is a chain in $P$ with $p_{i,1}\lessdot\cdots\lessdot p_{i,n_i}$;
 	\item
 	$p_{1,1}$ and $p_{2,1}$ are incomparable in $P$;
 	\item
 	$p_{1,n_1}$ and $p_{2,n_2}$ are incomparable in $P$.
 \end{enumerate}
 Let $P'$ be the induced subposet of $P$ on the set $\mathcal{S}$. If $R[\MI(P')]$ does not satisfy property $N_2$ then so  does $R[\MI(P)]$.
 \end{Theorem}
 
\begin{proof}
For $P$, let $P_1,\ldots,P_m,P_{m+1},\ldots,P_s$ be as defined in Discussion~\ref{reducesubposet}. For $1\leq i \leq s$, let $P'_i$ be the induced subposet of $P_i$ on the set $\mathcal{S}$. For  $1\leq i\leq s-1$, every order relation between the elements of $\mathcal{S}$ in $P_i$ is also an order relation in $P_{i+1}$. Also, $\{p_{1,1},p_{2,1}\}$ and $\{p_{1,n_{1}},p_{2,n_2}\}$ are the sets of minimal and maximal elements of $P_i$ respectively, for all $i = 1,\ldots, s$. Therefore, by Lemma~\ref{addingrelations}, $\beta_{24}(R[\MI(P_i')]) \neq 0$ for all $1\leq i \leq s$. By Discussion~\ref{reducesubposet}, it is enough to show that $R[\MI(P_s)]$ does not satisfy property $N_2$. We may replace $P$ by $P_s$ and $P'$ by $P_s'$.\par

Let $P^\partial$ be the dual poset of $P$. If $q\in P\setminus \mathcal{S}$, then   $\#\mathcal{S}^P_{q} \geq 2$. So if there exists a $q\in P^\partial\setminus \mathcal{S}$ with $\# \mathcal{S}^{P^\partial}_q \geq 2$, then $P$ contains a cover-preserving subposet as shown in Figure~\ref{fig5}. Thus, by Lemma~\ref{rank3simple}, $R[\MI(P)]$ does not satisfy property $N_2$.  So we may assume that for all $p\in P^\partial\setminus \mathcal{S}$, $\# \mathcal{S}^{P^\partial}_p=1$. Repeating the argument of Discussion~\ref{reducesubposet} for $P^{\partial}$, we obtain a poset $Q$ such that there is no $q\in Q\setminus \mathcal{S}$ with $\#\mathcal{S}^{Q	}_q = 1$. Observe that $Q$ is a  poset on the set $\mathcal{S}$. By Discussion~\ref{reducesubposet}, it suffices to prove that $R[\MI(Q)]$ does not satisfy property $N_2$. Note that $Q^\partial$ is a poset on the set $\mathcal{S}$ and all order relations of $P'$ are  also the order relations of $Q^\partial$. So by Lemma~\ref{addingrelations}, $R[\MI(Q^\partial)]$ does not satisfy property $N_2$. Thus, by Lemma~\ref{reflection}, $R[\MI(Q)]$ does not satisfy property $N_2$. Hence the proof.
\end{proof}

\begin{Remark}\label{thmplanar}
Note that in the proof of Theorem~\ref{planar}, the reduction from the poset $P$ to the poset $Q^\partial$ is independent of the hypothesis that $\MI(P')$ is a planar distributive lattice. In fact, we will also use the reduction from $P$ to $Q^\partial$ in Discussion~\ref{reduction3chains} where the distributive lattice is not restricted to be planar.

We have only used the assumption $\MI(P')$ is a planar distributive lattice to conclude that $\beta_{24}(R[\MI(Q^\partial)]) \neq 0$ and $\beta_{24}(R[\MI(P_i')]) \neq 0$ for all $i= 1,\ldots,s$. 
\end{Remark}

\subsection{}  
In this subsection, we prove a result analogous to Ene's result. Suppose that a poset can be decomposed into a union of three chains and it has three maximal and minimal elements. We prove some necessary conditions regarding when the Hibi ring associated to such poset satisfies property $N_2$.

\begin{Lemma}\label{poset9elements}
   Let  P be a poset on the disjoint union $\mathop{\cup}_{i=1}^{3} \{p_{i,1},p_{i,2}, p_{i,3}\}$ such that
   \begin{enumerate}
   	\item 
   	for all $1\leq i \leq 3$, $\{p_{i,1},p_{i,2},p_{i,3} \}$  is a chain in $P$ with $p_{i,1}\lessdot p_{i,2}\lessdot p_{i,3}$;
   	\item
   	$\{p_{1,1},p_{2,1},p_{3,1}\}$ and $\{p_{1,3},p_{2,3},p_{3,3}\}$ are the sets of  minimal  and maximal elements of $P$ respectively.
   \end{enumerate}
If $P$ is pure and connected, then $\beta_{24}(R[\MI(P)]) \neq 0$.
\end{Lemma}

Now we prove the main theorem of this subsection.

\begin{Theorem}\label{connected}
	Let  P be a poset on the set $\mathop{\cup}_{i=1}^{3} \{p_{i,1},..., p_{i,n_i}\}$ such that
	\begin{enumerate}
		\item 
		$p_{1,1},p_{2,1},p_{3,1}$ are distinct and $p_{1,n_{1}},p_{2,n_{2}},p_{3,n_{3}}$ are distinct;
		\item
		$\{p_{1,1},p_{2,1},p_{3,1}\}$ and $\{p_{1,n_{1}},p_{2,n_{2}},p_{3,n_{3}}\}$ are the sets of  minimal  and maximal elements of $P$ respectively;
		\item
		for all $1\leq i \leq 3$, $n_i\geq 3$; $\{p_{i,1},\ldots,p_{i,n_{i}} \}$  is a chain in $P$ with $p_{i,1}\lessdot \cdots\lessdot p_{i,n_{i}}$.		
	\end{enumerate} 
	If $P$ is connected and none of the minimal elements of $P$ is covered by a maximal element, then $\beta_{24}(R[\MI(P)]) \neq 0$.
\end{Theorem}

\begin{proof}
	Reduce $P$ to $P_n$,  where $P_n$ is as defined in Discussion~\ref{reducingposet}.  Since $P$ is connected, so is $P_n$. Since none of the minimal elements of $P$ is covered by a maximal element, we obtain that $P_n$ is pure. So by Discussion~\ref{reducingposet}, we may replace $P$ by $P_n$ and assume that $P$ is pure and $n_i = 3$ for all $1\leq i \leq 3$. Let $X'$ be as defined in Notation~\ref{construction1}. We will prove the result in the following cases:
	\begin{enumerate}
		\item[$(1)$] If $\# X'= 1$, then the result follows from Lemma~\ref{nonsimple}.
		\item[$(2)$] If $\# X'= 2$, then $P$ will contain  a cover-preserving subposet as shown in Figure~\ref{fig5}. Hence, the result follows from Lemma~\ref{rank3simple}. 
		\item[$(3)$] If $\# X'= 3$, then the result follows from Lemma~\ref{poset9elements}.
	\end{enumerate}  
\end{proof} 

Now, we proceed to prove Lemma~\ref{poset9elements}. Before that we prove some relevant lemmas. 

\begin{Lemma}\label{covering3elements}
	Let  P be as defined in Lemma~\ref{poset9elements}. If there exists an element in $P$ such that either it cover three elements or it is covered by three elements and $P$ is not as shown in figure~\ref{fig30}, then $\beta_{24}(R[\MI(P)]) \neq 0$.  
\end{Lemma}

\begin{proof}
By Theorem~\ref{planar}, we may assume that there is no subposet $P'$ of $P$, as defined in Theorem~\ref{planar}, with $\beta_{24}(R[\MI(P')]) \neq 0$.
Let $p\in P$ be the element such that  either it cover three elements or it is covered by three elements. Note that $p$ is either a maximal element or a minimal element or $p\in \{p_{1,2},p_{2,2},p_{3,2}\}$. If $p$ is a maximal element of $P$, then in $P^{\partial}$, $p$ is a minimal element and it is covered by three elements. In this case by Lemma~\ref{reflection}, replace $P$ by $P^\partial$ and we may assume that $p$ is a minimal element of $P$. So we only have to consider the cases when either $p$ is a minimal element or $p\in \{p_{1,2},p_{2,2},p_{3,2}\}$.  In most of the subcases of these two cases, we show that there exist $\delta, \gamma\in \MI(P)$ such that $\beta_{24}(R[L']) \neq 0$, where $L'$ is the sublattice $\{\alpha\in \MI(P): \delta\leq\alpha\leq\gamma \}$. Hence, by Proposition~\ref{semigroup} and Proposition~\ref{homologicallypure}, we conclude that $\beta_{24}(R[\MI(P)]) \neq 0$.
\begin{itemize}
  \item[$ Case\ 1$] Assume that $p$ is a minimal element of $P$. Possibly by  relabelling the elements of $P$, we may assume that $p=p_{1,1}$. We consider the following two subcases:\par
 
 \item[$(a)$] Consider the subcase when $p_{1,1}$ is covered by $\{p_{1,2},p_{2,2},p_{3,2}\}$. 
 If either $p_{3,3}$ covers $p_{2,2}$ or $p_{1,3}$ covers any element other than $p_{1,2}$, then exists a subposet $P'$ of $P$, as defined in Theorem~\ref{planar}, with $\beta_{24}(R[\MI(P')]) \neq 0$. So we may assume that $p_{3,3}$ does not cover $p_{2,2}$ and $p_{1,3}$ covers $p_{1,2}$ only. We proceed in the following two subsubcases:\\
 $(i)$ Assume that $p_{3,3}$ is only covering $p_{1,2}$ and $p_{3,2}$. Observe that $\delta=\emptyset$ and $\gamma = P\setminus\{p_{2,3}\}$ are the order ideals of $P$. By Proposition~\ref{sublattice}, $L' \cong \MI(P')$, where $P'$ is the poset as shown in Figure~\ref{fig26}.  One can use a computer to check that  $\beta_{24}(R[\MI(P')]) \neq 0$.\\
 $(ii)$ Now, assume that either $p_{3,3}$ is covering $p_{3,2}$ only or $p_{3,3}$ is  covering atleast $p_{2,1}$ and $p_{3,2}$. Let $\delta=\emptyset$ and $\gamma = P\setminus\{p_{1,3},p_{2,3}\}$. By Proposition~\ref{sublattice}, $L' \cong \MI(P')$, where $P'$ is one of the posets as shown in Figure~\ref{fig25}-\ref{fig28}.  Again, it can be checked by a computer that  $\beta_{24}(R[\MI(P')]) \neq 0$ for all $P'$.\par 
 
 \item[$(b)$] Consider the subcase when $p_{1,1}$ is not covered by $\{p_{1,2},p_{2,2},p_{3,2}\}$. So $p_{1,1}$ is either covered by  $\{p_{1,2},p_{2,3},p_{3,2}\}$ or $\{p_{1,2},p_{2,2},p_{3,3}\}$ or $\{p_{1,2},p_{2,3},p_{3,3}\}$. By symmetry, it is enough to consider one of the cases from $\{p_{1,2},p_{2,3},p_{3,2}\}$ and $\{p_{1,2},p_{2,2},p_{3,3}\}$.\par 
 
 First, consider the subsubcase when $p_{1,1}$ is covered by $\{p_{1,2},p_{2,3},p_{3,2}\}$. We have $p_{2,1}\lessdot p_{2,2}$, reduce $P$ to $P_1$ using the methods of  Discussion~\ref{construction1}. If $p_{2,1}$ is covered by $p_{1,2}$ or $p_{3,2}$ in $P$, then $P_1$ will contain a cover-preserving subposet as shown in Figure~\ref{fig4}. So we may assume that $p_{2,1}$ is not covered by $p_{1,2}$ and $p_{3,2}$. Observe that $P_1$ is a poset on the underlying set $P\setminus \{p_{2,1}\}$. Also, $\{p_{1,1},p_{2,2},p_{3,1}\}$ and $\{p_{1,3},p_{2,3},p_{3,3}\}$ are the sets of minimal and maximal elements of $P_1$ respectively. Also, $p_{1,1}$ is covered by $\{p_{1,2},p_{2,3},p_{3,2}\}$ in $P_1$. Repeating the argument of the subcase $(a)$, we deduce that the result holds in this subsubcase.\par
 
 Now, we consider the subsubcase when $p_{1,1}$ is covered by $\{p_{1,2},p_{2,3},p_{3,3}\}$. Again, we have $p_{2,1}\lessdot p_{2,2}$, reduce $P$ to $P_1$ using the methods of  Discussion~\ref{construction1}. If $p_{2,1}$ is covered by $p_{1,2}$  in $P$, then $P_1$ will contain a cover-preserving subposet as shown in Figure~\ref{fig4}. So we may assume that $p_{2,1}$ is not covered by $p_{1,2}$ in $P$. Similarly, we may assume that $p_{3,1}$ is not covered by $p_{1,2}$. If either $p_{2,2}$ or $p_{3,2}$ is covered by $p_{1,3}$, then $P$ will contain a subposet $P'$, as defined in Theorem~\ref{planar},  with $\beta_{24}(R[\MI(P')]) \neq 0$. If either $p_{2,2}$ is covered by $p_{3,3}$ or $p_{3,2}$ is covered by $p_{2,3}$, then we are done by the previous subsubcase. Since $P$ is not as shown in figure~\ref{fig30}, the only possibility for $P$ is that $P$ is isomorphic to one of the posets as shown in Figure~\ref{fig31}-\ref{fig37}.  One can use a computer to check that  $\beta_{24}(R[\MI(P)]) \neq 0$.

 \item[$ Case\ 2$] Assume that $p\in \{p_{1,2},p_{2,2},p_{3,2}\}$. Possibly by replacing $P$ with $P^\partial$, we may assume that $p$ is covering all the minimal  elements. Possibly by  relabelling the elements of $P$, we may assume that $p=p_{1,2}$. If $p_{2,2}$ and $p_{3,2}$ are covered by $p_{1,3}$, then we are done by Case 1. So we may assume that both $p_{2,2}$ and $p_{3,2}$ are not covered by $p_{1,3}$. Let $\delta=\emptyset$ and $\gamma = P\setminus\{p_{2,3},p_{3,3}\}$. By Proposition~\ref{sublattice}, $L' \cong \MI(P')$, where $P'$ is one of the posets as shown in Figure~\ref{fig23}-\ref{fig24}.  Again, one can use a computer to check that  $\beta_{24}(R[\MI(P')]) \neq 0$.
\end{itemize} 
\end{proof}
 	
\begin{figure}[h]
	\begin{subfigure}[t]{4cm}
		\centering	
		\begin{tikzpicture}[scale=2]
		\draw[] (0,0)--(0,.6)--(0,1.2) (0,0)--(.4,1.2)--(.4,.6)--(.4,0) 
		(.7,0)--(.7,.6)--(.7,1.2)--(0,0) ;
		\filldraw[black] (0,0) circle (.5pt) node[anchor=north]  {$p_1$}; 
		\filldraw[black] (.7,0) circle (.5pt) node[anchor=north] {$p_3$}; 
		\filldraw[black] (.4,0) circle (.5pt) node[anchor=north] {$p_2$}; 
		\filldraw[black] (0,.6) circle (.5pt) node[anchor=east]  {$p_4$};   
		\filldraw[black] (.4,.5) circle (.5pt) node[anchor=west]  {$p_5$}; 		
		\filldraw[black] (.7,.6) circle (.5pt) node[anchor=west]  {$p_{6}$}; 
		\filldraw[black] (.7,1.2) circle (.5pt) node[anchor=south] {$p_{9}$};   
		\filldraw[black] (.0,1.2) circle (.5pt) node[anchor=south] {$p_{7}$}; 		\filldraw[black] (.4,1.2) circle (.5pt) node[anchor=south]{$p_{8}$};   
		\end{tikzpicture}
		\caption{}\label{fig30}
	\end{subfigure} 
	\quad	
    \begin{subfigure}[t]{4cm}
		\centering	
		\begin{tikzpicture}[scale=2]
		\draw[] (0,0)--(0,.6)--(0,1.2) (0,0.6)--(.4,0)--(.4,.6) 
		(0,.6)--(.8,0)--(.8,.6);
		\filldraw[black] (0,0) circle (.5pt) node[anchor=north]  {$p_1$}; 
		\filldraw[black] (.8,0) circle (.5pt) node[anchor=north] {$p_3$}; 
		\filldraw[black] (.4,0) circle (.5pt) node[anchor=north] {$p_2$}; 
		\filldraw[black] (0,.6) circle (.5pt) node[anchor=east]  {$p_4$};   
		\filldraw[black] (.4,.6) circle (.5pt) node[anchor=south]  {$p_5$}; 		
		\filldraw[black] (.8,.6) circle (.5pt) node[anchor=south]  {$p_{6}$}; 
		\filldraw[black] (0,1.2) circle (.5pt) node[anchor=south] {$p_{7}$};   
		
		\end{tikzpicture}
		\caption{$\beta_{24}=8$}\label{fig23}
	\end{subfigure}  	
	\quad	
	\begin{subfigure}[t]{4cm}
		\centering	
		\begin{tikzpicture}[scale=2]
		\draw[] (0,0)--(0,.6)--(0,1.2)--(.4,.6) (0,0.6)--(.4,0)--(.4,.6) 
		(0,.6)--(.8,0)--(.8,.6);
		\filldraw[black] (0,0) circle (.5pt) node[anchor=north]  {$p_1$}; 
		\filldraw[black] (.8,0) circle (.5pt) node[anchor=north] {$p_3$}; 
		\filldraw[black] (.4,0) circle (.5pt) node[anchor=north] {$p_2$}; 
		\filldraw[black] (0,.6) circle (.5pt) node[anchor=east]  {$p_4$};   
		\filldraw[black] (.4,.6) circle (.5pt) node[anchor=south]  {$p_5$}; 		
		\filldraw[black] (.8,.6) circle (.5pt) node[anchor=south]  {$p_{6}$}; 
		\filldraw[black] (0,1.2) circle (.5pt) node[anchor=south] {$p_{7}$};  
		
		\end{tikzpicture}
		\caption{$\beta_{24}=42$}\label{fig24}
	\end{subfigure}  	 	
	\quad	
	\begin{subfigure}[t]{4cm}
		\centering	
		\begin{tikzpicture}[scale=2]
		\draw[] (0,0)--(0,.6)--(0,1.2) (0,0)--(.4,.6)--(.4,0) 
		(.8,0)--(.8,.6)--(0,0) (.8,.6)--(.8,1.2)--(0,.6);
		\filldraw[black] (0,0) circle (.5pt) node[anchor=north]  {$p_1$}; 
		\filldraw[black] (.8,0) circle (.5pt) node[anchor=north] {$p_3$}; 
		\filldraw[black] (.4,0) circle (.5pt) node[anchor=north] {$p_2$}; 
		\filldraw[black] (0,.6) circle (.5pt) node[anchor=east]  {$p_4$};   
		\filldraw[black] (.4,.6) circle (.5pt) node[anchor=south]  {$p_5$}; 		
		\filldraw[black] (.8,.6) circle (.5pt) node[anchor=west]  {$p_{6}$}; 
		\filldraw[black] (.8,1.2) circle (.5pt) node[anchor=south] {$p_{8}$};  
		\filldraw[black] (.0,1.2) circle (.5pt) node[anchor=south] {$p_{7}$};  
		\end{tikzpicture}
		\caption{$\beta_{24}=1$}\label{fig26}
	\end{subfigure} 
	\quad		
	\begin{subfigure}[t]{4cm}
		\centering	
		\begin{tikzpicture}[scale=2]
		\draw[] (0,0)--(0,.6) (0,0)--(.4,.6)--(.4,0) 
		(.8,0)--(.8,.6)--(0,0) (.8,.6)--(.8,1.2);
		\filldraw[black] (0,0) circle (.5pt) node[anchor=north]  {$p_1$}; 
		\filldraw[black] (.8,0) circle (.5pt) node[anchor=north] {$p_3$}; 
		\filldraw[black] (.4,0) circle (.5pt) node[anchor=north] {$p_2$}; 
		\filldraw[black] (0,.6) circle (.5pt) node[anchor=south]  {$p_4$};   
		\filldraw[black] (.4,.6) circle (.5pt) node[anchor=south]  {$p_5$}; 		
		\filldraw[black] (.8,.6) circle (.5pt) node[anchor=west]  {$p_{6}$}; 
		\filldraw[black] (.8,1.2) circle (.5pt) node[anchor=south] {$p_{7}$};   
		\end{tikzpicture}
		\caption{$\beta_{24}=1$}\label{fig25}
	\end{subfigure} 	 
   	\quad	
	\begin{subfigure}[t]{4cm}
		\centering	
		\begin{tikzpicture}[scale=2]
		\draw[] (0,0)--(0,.6) (0,0)--(.4,.6)--(.4,0) 
		(.8,0)--(.8,.6)--(0,0) (.8,.6)--(.8,1.2)--(0.4,0);
		\filldraw[black] (0,0) circle (.5pt) node[anchor=north]  {$p_1$}; 
		\filldraw[black] (.8,0) circle (.5pt) node[anchor=north] {$p_3$}; 
		\filldraw[black] (.4,0) circle (.5pt) node[anchor=north] {$p_2$}; 
		\filldraw[black] (0,.6) circle (.5pt) node[anchor=south]  {$p_4$};   
		\filldraw[black] (.4,.6) circle (.5pt) node[anchor=south]  {$p_5$}; 		
		\filldraw[black] (.8,.6) circle (.5pt) node[anchor=west]  {$p_{6}$}; 
		\filldraw[black] (.8,1.2) circle (.5pt) node[anchor=south] {$p_{7}$};   
		\end{tikzpicture}
		\caption{$\beta_{24}=42$}\label{fig27}
	\end{subfigure} 
	\quad	
	\begin{subfigure}[t]{4cm}
		\centering	
		\begin{tikzpicture}[scale=2]
		\draw[] (0,0)--(0,.6) (0,0)--(.4,.6)--(.4,0) 
		(.8,0)--(.8,.6)--(0,0) (.8,.6)--(.8,1.2)--(0.4,0) (.8,1.2)--(0,.6);
		\filldraw[black] (0,0) circle (.5pt) node[anchor=north]  {$p_1$}; 
		\filldraw[black] (.8,0) circle (.5pt) node[anchor=north] {$p_3$}; 
		\filldraw[black] (.4,0) circle (.5pt) node[anchor=north] {$p_2$}; 
		\filldraw[black] (0,.6) circle (.5pt) node[anchor=south]  {$p_4$};   
		\filldraw[black] (.4,.6) circle (.5pt) node[anchor=south]  {$p_5$}; 		
		\filldraw[black] (.8,.6) circle (.5pt) node[anchor=west]  {$p_{6}$}; 
		\filldraw[black] (.8,1.2) circle (.5pt) node[anchor=south] {$p_{7}$};   
		\end{tikzpicture}
		\caption{$\beta_{24}=28$}\label{fig28}
	\end{subfigure} 
	\quad		
	\begin{subfigure}[t]{4cm}
		\centering	
		\begin{tikzpicture}[scale=2]
		\draw[] (0,0)--(0,.6)--(0,1.2) (0,0)--(.4,1.2)--(.4,.6)--(.4,0) 
		(0.4,.5)--(.7,0)--(.7,.6)--(.7,1.2)--(0,0) ;
		\filldraw[black] (0,0) circle (.5pt) node[anchor=north]  {$p_1$}; 
		\filldraw[black] (.7,0) circle (.5pt) node[anchor=north] {$p_3$}; 
		\filldraw[black] (.4,0) circle (.5pt) node[anchor=north] {$p_2$}; 
		\filldraw[black] (0,.6) circle (.5pt) node[anchor=east]  {$p_4$};   
		\filldraw[black] (.4,.5) circle (.5pt) node[anchor=west]  {$p_5$}; 		
		\filldraw[black] (.7,.6) circle (.5pt) node[anchor=west]  {$p_{6}$}; 
		\filldraw[black] (.7,1.2) circle (.5pt) node[anchor=south] {$p_{9}$};   
		\filldraw[black] (.0,1.2) circle (.5pt) node[anchor=south] {$p_{7}$}; 		\filldraw[black] (.4,1.2) circle (.5pt) node[anchor=south]{$p_{8}$}; 				
		\end{tikzpicture}
		\caption{$\beta_{24}=299$}\label{fig31}
	\end{subfigure} 
	\quad
	\begin{subfigure}[t]{4cm}
		\centering	
		\begin{tikzpicture}[scale=2]
		\draw[] (0,0)--(0,.6)--(0,1.2) (0,0)--(.4,1.2)--(.4,.6)--(.4,0) 
		(0.4,1.2)--(.7,0)--(.7,.6)--(.7,1.2)--(0,0) ;
		\filldraw[black] (0,0) circle (.5pt) node[anchor=north]  {$p_1$}; 
		\filldraw[black] (.7,0) circle (.5pt) node[anchor=north] {$p_3$}; 
		\filldraw[black] (.4,0) circle (.5pt) node[anchor=north] {$p_2$}; 
		\filldraw[black] (0,.6) circle (.5pt) node[anchor=east]  {$p_4$};   
		\filldraw[black] (.4,.5) circle (.5pt) node[anchor=west]  {$p_5$}; 		
		\filldraw[black] (.7,.6) circle (.5pt) node[anchor=west]  {$p_{6}$}; 
		\filldraw[black] (.7,1.2) circle (.5pt) node[anchor=south] {$p_{9}$};   
		\filldraw[black] (.0,1.2) circle (.5pt) node[anchor=south] {$p_{7}$}; 		\filldraw[black] (.4,1.2) circle (.5pt) node[anchor=south]{$p_{8}$}; 				
		\end{tikzpicture}
		\caption{$\beta_{24}=327$}\label{fig37}
	\end{subfigure}				  		     	 	   	
	\caption{}\label{fig22}
\end{figure}

\begin{Lemma}\label{subposetconnected}
Let  P be as defined in Lemma~\ref{poset9elements}. If the induced subposet of $P$, defined on the underlying set $P\setminus \{p_{1,1},p_{2,1}, p_{3,1}\}$ or $P\setminus \{p_{1,3},p_{2,3}, p_{3,3}\}$, is connected. Then $\beta_{24}(R[\MI(P)]) \neq 0$.
\end{Lemma}
\begin{proof}
By Theorem~\ref{planar}, we may assume that there is no subposet $P'$ of $P$, as defined in Theorem~\ref{planar}, with $\beta_{24}(R[\MI(P')]) \neq 0$. Possibly by replacing $P$ with $P^\partial$, we may assume that the subposet $P'$ of $P$ defined on the underlying set $P\setminus \{p_{1,3},p_{2,3}, p_{3,3}\}$ is connected. Observe that $P'$ is isomorphic to one of the posets as shown in Figure~\ref{fig33}-\ref{fig34}. If $P'$ is as shown in Figure~\ref{fig34}, then we are done by Lemma~\ref{covering3elements}. 

Now, consider the case when $P'$ is as shown in Figure~\ref{fig33}. Possibly by  relabelling the elements of $P$, we may assume that $p_{1,2}$ is covering exactly one minimal element of $P$. If either $p_{1,1}\lessdot p_{3,3}$ or $p_{3,1}\lessdot p_{2,3}$ or $p_{2,2}\lessdot p_{3,3}$ or $p_{3,2}\lessdot p_{2,3}$, then there exists a subposet $P'$ of $P$, as defined in Theorem~\ref{planar}, with $\beta_{24}(R[\MI(P')]) \neq 0$. Let $\delta=\emptyset$ and $\gamma = P\setminus\{p_{1,2},p_{1,3}\}$. Let also $L' = \{\alpha\in \MI(P): \delta\leq\alpha\leq\gamma \}$. By Proposition~\ref{sublattice}, $L' \cong \MI(P_1)$, where $P_1$ is as shown in Figure~\ref{fig29}. One can use a computer to check that  $\beta_{24}(R[\MI(P_1)]) \neq 0$. Hence, by Proposition~\ref{semigroup} and Proposition~\ref{homologicallypure}, $\beta_{24}(R[\MI(P)]) \neq 0$.
\end{proof}
\begin{figure}
	\begin{subfigure}[t]{4cm}
		\centering	
		\begin{tikzpicture}[scale=2]
		\draw[] (0,0)--(0,.6) (0,0)--(.4,.6)--(.4,0)--(.8,.6)	(.8,0)--(.8,.6) ;
		\filldraw[black] (0,0) circle (.5pt) node[anchor=north]  {$p_1$}; 
		\filldraw[black] (.8,0) circle (.5pt) node[anchor=north] {$p_3$}; 
		\filldraw[black] (.4,0) circle (.5pt) node[anchor=north] {$p_2$}; 
		\filldraw[black] (0,.6) circle (.5pt) node[anchor=south]  {$p_4$};   
		\filldraw[black] (.4,.6) circle (.5pt) node[anchor=south]  {$p_5$}; 		
		\filldraw[black] (.8,.6) circle (.5pt) node[anchor=south]  {$p_{6}$};
				
		\end{tikzpicture}
		\caption{}\label{fig33}
	\end{subfigure}
	\quad	
	\begin{subfigure}[t]{4cm}
		\centering	
		\begin{tikzpicture}[scale=2]
		\draw[] (0,0)--(0,.6) (0,0)--(.4,.6)--(.4,0)	(.8,0)--(.8,.6)--(0,0) ;
		\filldraw[black] (0,0) circle (.5pt) node[anchor=north]  {$p_1$}; 
		\filldraw[black] (.8,0) circle (.5pt) node[anchor=north] {$p_3$}; 
		\filldraw[black] (.4,0) circle (.5pt) node[anchor=north] {$p_2$}; 
		\filldraw[black] (0,.6) circle (.5pt) node[anchor=south]  {$p_4$};   
		\filldraw[black] (.4,.6) circle (.5pt) node[anchor=south]  {$p_5$}; 		
		\filldraw[black] (.8,.6) circle (.5pt) node[anchor=south]  {$p_{6}$};  				
		\end{tikzpicture}
		\caption{}\label{fig34}
	\end{subfigure}		
	\quad
	\begin{subfigure}[t]{4cm}
		\centering	
		\begin{tikzpicture}[scale=2]
		\draw[] (.8,0)--(.4,.6)--(0,0) (.4,.6)--(.4,1.1) 
		(.8,0)--(1,.6)--(1,1.1)  (1.4,0)--(1,.6);
		\filldraw[black] (0,0) circle (.5pt) node[anchor=north]  {$p_1$}; 
		\filldraw[black] (1.4,0) circle (.5pt) node[anchor=north] {$p_3$}; 
		\filldraw[black] (.8,0) circle (.5pt) node[anchor=north] {$p_2$}; 
		\filldraw[black] (.4,.6) circle (.5pt) node[anchor=east]  {$p_4$};   
		\filldraw[black] (1,.6) circle (.5pt) node[anchor=east]  {$p_5$}; 		
		\filldraw[black] (.4,1.1) circle (.5pt) node[anchor=south]  {$p_{6}$}; 
		\filldraw[black] (1,1.1) circle (.5pt) node[anchor=south] {$p_{7}$};   
		\end{tikzpicture}
		\caption{$\beta_{24}=1$}\label{fig29}
	\end{subfigure}
	\quad		
	\caption{}\label{fig32}
\end{figure}

\begin{proof}[Proof of Lemma~\ref{poset9elements}]
By Theorem~\ref{planar}, we may assume that there is no subposet $P'$ of $P$, as defined in Theorem~\ref{planar}, with $\beta_{24}(R[\MI(P')]) \neq 0$. By Lemma~\ref{covering3elements}, we may assume that there is no element in $P$ such that either it cover three elements or it is covered by three elements. By Lemma~\ref{subposetconnected}, we may assume that the subposets of $P$ defined on the underlying sets $P\setminus \{p_{1,1},p_{2,1}, p_{3,1}\}$ and $P\setminus \{p_{1,3},p_{2,3}, p_{3,3}\}$ are not connected. Then $P$ is isomorphic to one of the posets as shown in Figure~\ref{fig35}. One can use a computer to check that  $\beta_{24}(R[\MI(P)]) \neq 0$. This concludes the proof.
\end{proof}
\begin{figure}
	
	\begin{subfigure}[t]{4cm}
		\centering	
		\begin{tikzpicture}[scale=2]
		\draw[] (0,0)--(0,.6)--(0,1.2) (.4,1.2)--(.4,.6)--(.4,0) (.4,.6)--(0,0)
		(.8,0)--(.8,.6)--(.8,1.2)--(0,.6) ;
		\filldraw[black] (0,0) circle (.5pt) node[anchor=north]  {$p_1$}; 
		\filldraw[black] (.8,0) circle (.5pt) node[anchor=north] {$p_3$}; 
		\filldraw[black] (.4,0) circle (.5pt) node[anchor=north] {$p_2$}; 
		\filldraw[black] (0,.6) circle (.5pt) node[anchor=east]  {$p_4$};   
		\filldraw[black] (.4,.6) circle (.5pt) node[anchor=west]  {$p_5$}; 		
		\filldraw[black] (.8,.6) circle (.5pt) node[anchor=west]  {$p_{6}$}; 
		\filldraw[black] (.8,1.2) circle (.5pt) node[anchor=south] {$p_{9}$};   
		\filldraw[black] (.0,1.2) circle (.5pt) node[anchor=south] {$p_{7}$}; 		\filldraw[black] (.4,1.2) circle (.5pt) node[anchor=south]{$p_{8}$}; 				
		\end{tikzpicture}
		\caption{$\beta_{24}=3$}\ 
	\end{subfigure} 	
	\quad
	\begin{subfigure}[t]{4cm}
		\centering	
		\begin{tikzpicture}[scale=2]
		\draw[] (0,0)--(0,.6)--(0,1.2) (.4,1.2)--(.4,.6)--(.4,0) (.4,.6)--(.8,0)
		(.8,0)--(.8,.6)--(.8,1.2)--(0,.6) ;
		\filldraw[black] (0,0) circle (.5pt) node[anchor=north]  {$p_1$}; 
		\filldraw[black] (.8,0) circle (.5pt) node[anchor=north] {$p_3$}; 
		\filldraw[black] (.4,0) circle (.5pt) node[anchor=north] {$p_2$}; 
		\filldraw[black] (0,.6) circle (.5pt) node[anchor=east]  {$p_4$};   
		\filldraw[black] (.4,.6) circle (.5pt) node[anchor=west]  {$p_5$}; 		
		\filldraw[black] (.8,.6) circle (.5pt) node[anchor=west]  {$p_{6}$}; 
		\filldraw[black] (.8,1.2) circle (.5pt) node[anchor=south] {$p_{9}$};   
		\filldraw[black] (.0,1.2) circle (.5pt) node[anchor=south] {$p_{7}$}; 		\filldraw[black] (.4,1.2) circle (.5pt) node[anchor=south]{$p_{8}$}; 				
		\end{tikzpicture}
		\caption{$\beta_{24}=1$}\ 
	\end{subfigure} 
	\quad
	\begin{subfigure}[t]{4cm}
		\centering	
		\begin{tikzpicture}[scale=2]
		\draw[] (0,0)--(0,.6)--(0,1.2) (.8,.6)--(.4,1.2)--(.4,.6)--(.4,0) (.4,.6)--(0,0)	(.8,0)--(.8,.6)--(.8,1.2) ;
		\filldraw[black] (0,0) circle (.5pt) node[anchor=north]  {$p_1$}; 
		\filldraw[black] (.8,0) circle (.5pt) node[anchor=north] {$p_3$}; 
		\filldraw[black] (.4,0) circle (.5pt) node[anchor=north] {$p_2$}; 
		\filldraw[black] (0,.6) circle (.5pt) node[anchor=east]  {$p_4$};   
		\filldraw[black] (.4,.6) circle (.5pt) node[anchor=west]  {$p_5$}; 		
		\filldraw[black] (.8,.6) circle (.5pt) node[anchor=west]  {$p_{6}$}; 
		\filldraw[black] (.8,1.2) circle (.5pt) node[anchor=south] {$p_{9}$};   
		\filldraw[black] (.0,1.2) circle (.5pt) node[anchor=south] {$p_{7}$}; 		\filldraw[black] (.4,1.2) circle (.5pt) node[anchor=south]{$p_{8}$}; 				
		\end{tikzpicture}
		\caption{$\beta_{24}=3$}\ 
	\end{subfigure}
	\caption{}\label{fig35}
\end{figure}

\begin{figure}
	
	\begin{subfigure}[t]{4cm}
		\centering	
		\begin{tikzpicture}[scale=2]
		\draw[] (0,0)--(0,.6)--(0,1.2) (.8,0)--(.4,1.2)--(.4,.6)--(.4,0) (.4,.6)--(0,0)	(.8,0)--(.8,.6)--(.8,1.2) ;
		\filldraw[black] (0,0) circle (.5pt) node[anchor=north]  {$p_1$}; 
		\filldraw[black] (.8,0) circle (.5pt) node[anchor=north] {$p_3$}; 
		\filldraw[black] (.4,0) circle (.5pt) node[anchor=north] {$p_2$}; 
		\filldraw[black] (0,.6) circle (.5pt) node[anchor=east]  {$p_4$};   
		\filldraw[black] (.4,.6) circle (.5pt) node[anchor=east]  {$p_5$}; 		
		\filldraw[black] (.8,.6) circle (.5pt) node[anchor=west]  {$p_{6}$}; 
		\filldraw[black] (.8,1.2) circle (.5pt) node[anchor=south] {$p_{9}$};   
		\filldraw[black] (.0,1.2) circle (.5pt) node[anchor=south] {$p_{7}$}; 		\filldraw[black] (.4,1.2) circle (.5pt) node[anchor=south]{$p_{8}$}; 				
		\end{tikzpicture}
		\caption{$\beta_{24}=0$}\label{fig40}
	\end{subfigure}
	\quad
	\begin{subfigure}[t]{4cm}
		\centering	
		\begin{tikzpicture}[scale=2]
		\draw[] (0,0)--(0,.6)--(0,1.2) (0,0)--(.4,1.2)--(.4,.6)--(.4,0)--(.8,.6)	(.8,0)--(.8,.6)--(.8,1.2) ;
		\filldraw[black] (0,0) circle (.5pt) node[anchor=north]  {$p_1$}; 
		\filldraw[black] (.8,0) circle (.5pt) node[anchor=north] {$p_3$}; 
		\filldraw[black] (.4,0) circle (.5pt) node[anchor=north] {$p_2$}; 
		\filldraw[black] (0,.6) circle (.5pt) node[anchor=east]  {$p_4$};   
		\filldraw[black] (.4,.6) circle (.5pt) node[anchor=west]  {$p_5$}; 		
		\filldraw[black] (.8,.6) circle (.5pt) node[anchor=west]  {$p_{6}$}; 
		\filldraw[black] (.8,1.2) circle (.5pt) node[anchor=south] {$p_{9}$};   
		\filldraw[black] (.0,1.2) circle (.5pt) node[anchor=south] {$p_{7}$}; 		\filldraw[black] (.4,1.2) circle (.5pt) node[anchor=south]{$p_{8}$}; 				
		\end{tikzpicture}
		\caption{$\beta_{24}=0$}\ 
	\end{subfigure}	
	\quad	
	\begin{subfigure}[t]{4cm}
		\centering	
		\begin{tikzpicture}[scale=2]
		\draw[] (0,0)--(0,.6)--(0,1.2) (.4,1.2)--(.4,.6)--(.4,0) (.4,.6)--(0,0)	(.8,0)--(.8,.6)--(.8,1.2)--(.4,0) ;
		\filldraw[black] (0,0) circle (.5pt) node[anchor=north]  {$p_1$}; 
		\filldraw[black] (.8,0) circle (.5pt) node[anchor=north] {$p_3$}; 
		\filldraw[black] (.4,0) circle (.5pt) node[anchor=north] {$p_2$}; 
		\filldraw[black] (0,.6) circle (.5pt) node[anchor=east]  {$p_4$};   
		\filldraw[black] (.4,.6) circle (.5pt) node[anchor=east]  {$p_5$}; 		
		\filldraw[black] (.8,.6) circle (.5pt) node[anchor=west]  {$p_{6}$}; 
		\filldraw[black] (.8,1.2) circle (.5pt) node[anchor=south] {$p_{9}$};   
		\filldraw[black] (.0,1.2) circle (.5pt) node[anchor=south] {$p_{7}$}; 		\filldraw[black] (.4,1.2) circle (.5pt) node[anchor=south]{$p_{8}$}; 				
		\end{tikzpicture}
		\caption{$\beta_{24}=0$}\ 
	\end{subfigure}
	\quad
	\begin{subfigure}[t]{4cm}
		\centering	
		\begin{tikzpicture}[scale=2]
		\draw[] (0,0)--(0,.6)--(0,1.2) (0,0)--(.4,1.2)--(.4,.6)--(.4,0)	(.8,0)--(.8,.6)--(.8,1.2)--(.4,0) ;
		\filldraw[black] (0,0) circle (.5pt) node[anchor=north]  {$p_1$}; 
		\filldraw[black] (.8,0) circle (.5pt) node[anchor=north] {$p_3$}; 
		\filldraw[black] (.4,0) circle (.5pt) node[anchor=north] {$p_2$}; 
		\filldraw[black] (0,.6) circle (.5pt) node[anchor=east]  {$p_4$};   
		\filldraw[black] (.4,.6) circle (.5pt) node[anchor=east]  {$p_5$}; 		
		\filldraw[black] (.8,.6) circle (.5pt) node[anchor=west]  {$p_{6}$}; 
		\filldraw[black] (.8,1.2) circle (.5pt) node[anchor=south] {$p_{9}$};   
		\filldraw[black] (.0,1.2) circle (.5pt) node[anchor=south] {$p_{7}$}; 		\filldraw[black] (.4,1.2) circle (.5pt) node[anchor=south]{$p_{8}$}; 				
		\end{tikzpicture}
		\caption{$\beta_{24}=0$}\ 
	\end{subfigure}			
	\quad	
	\begin{subfigure}[t]{4cm}
		\centering	
		\begin{tikzpicture}[scale=2]
		\draw[] (.4,1.2)--(0,0)--(0,.6)--(0,1.2) 	(.4,1.2)--(.4,.6)--(.4,0)--(.8,.6)	(.8,0)--(.8,.6)--(.8,1.2)--(.4,.6)  ;
		\filldraw[black] (0,0) circle (.5pt) node[anchor=north]  {$p_1$}; 
		\filldraw[black] (.8,0) circle (.5pt) node[anchor=north] {$p_3$}; 
		\filldraw[black] (.4,0) circle (.5pt) node[anchor=north] {$p_2$}; 
		\filldraw[black] (0,.6) circle (.5pt) node[anchor=east]  {$p_4$};   
		\filldraw[black] (.4,.6) circle (.5pt) node[anchor=west]  {$p_5$}; 		
		\filldraw[black] (.8,.6) circle (.5pt) node[anchor=west]  {$p_{6}$}; 
		\filldraw[black] (.8,1.2) circle (.5pt) node[anchor=south] {$p_{9}$};   
		\filldraw[black] (.0,1.2) circle (.5pt) node[anchor=south] {$p_{7}$}; 		\filldraw[black] (.4,1.2) circle (.5pt) node[anchor=south]{$p_{8}$}; 				
		\end{tikzpicture}
		\caption{$\beta_{24}=0$}\label{fig41}
	\end{subfigure}
	\quad	
	\begin{subfigure}[t]{4cm}
		\centering	
		\begin{tikzpicture}[scale=2]
		\draw[] (0,0)--(0,.6)--(0,1.2) (.8,0)--(.4,1.2)--(.4,.6)--(.4,0) (.4,.6)--(0,0)
		(.8,0)--(.8,.6)--(.8,1.2)--(0,.6) ;
		\filldraw[black] (0,0) circle (.5pt) node[anchor=north]  {$p_1$}; 
		\filldraw[black] (.8,0) circle (.5pt) node[anchor=north] {$p_3$}; 
		\filldraw[black] (.4,0) circle (.5pt) node[anchor=north] {$p_2$}; 
		\filldraw[black] (0,.6) circle (.5pt) node[anchor=east]  {$p_4$};   
		\filldraw[black] (.4,.6) circle (.5pt) node[anchor=east]  {$p_5$}; 		
		\filldraw[black] (.8,.6) circle (.5pt) node[anchor=west]  {$p_{6}$}; 
		\filldraw[black] (.8,1.2) circle (.5pt) node[anchor=south] {$p_{9}$};   
		\filldraw[black] (.0,1.2) circle (.5pt) node[anchor=south] {$p_{7}$}; 		\filldraw[black] (.4,1.2) circle (.5pt) node[anchor=south]{$p_{8}$}; 				
		\end{tikzpicture}
		\caption{$\beta_{24}=285$}\label{fig42}
	\end{subfigure} 
	\quad	
	\begin{subfigure}[t]{4cm}
		\centering	
		\begin{tikzpicture}[scale=2]
		\draw[] (.4,.6)--(0,0)--(0,.6)--(0,1.2)--(.8,.6) 	(.4,1.2)--(.4,.6)--(.4,0)	(.4,1.2)--(.8,0)--(.8,.6)--(.8,1.2)  ;
		\filldraw[black] (0,0) circle (.5pt) node[anchor=north]  {$p_1$}; 
		\filldraw[black] (.8,0) circle (.5pt) node[anchor=north] {$p_3$}; 
		\filldraw[black] (.4,0) circle (.5pt) node[anchor=north] {$p_2$}; 
		\filldraw[black] (0,.6) circle (.5pt) node[anchor=east]  {$p_4$};   
		\filldraw[black] (.4,.6) circle (.5pt) node[anchor=east]  {$p_5$}; 		
		\filldraw[black] (.8,.6) circle (.5pt) node[anchor=west]  {$p_{6}$}; 
		\filldraw[black] (.8,1.2) circle (.5pt) node[anchor=south] {$p_{9}$};   
		\filldraw[black] (.0,1.2) circle (.5pt) node[anchor=south] {$p_{7}$}; 		\filldraw[black] (.4,1.2) circle (.5pt) node[anchor=south]{$p_{8}$}; 				
		\end{tikzpicture}
		\caption{$\beta_{24}=283$}\ 
	\end{subfigure}	
	\quad	
	\begin{subfigure}[t]{4cm}
		\centering	
		\begin{tikzpicture}[scale=2]
		\draw[] (.4,.6)--(0,0)--(0,.6)--(0,1.2)--(.8,.6) 	(.4,1.2)--(.4,.6)--(.4,0)	(.8,0)--(.8,.6)--(.8,1.2)--(.4,0)  ;
		\filldraw[black] (0,0) circle (.5pt) node[anchor=north]  {$p_1$}; 
		\filldraw[black] (.8,0) circle (.5pt) node[anchor=north] {$p_3$}; 
		\filldraw[black] (.4,0) circle (.5pt) node[anchor=north] {$p_2$}; 
		\filldraw[black] (0,.6) circle (.5pt) node[anchor=east]  {$p_4$};   
		\filldraw[black] (.4,.6) circle (.5pt) node[anchor=east]  {$p_5$}; 		
		\filldraw[black] (.8,.6) circle (.5pt) node[anchor=west]  {$p_{6}$}; 
		\filldraw[black] (.8,1.2) circle (.5pt) node[anchor=south] {$p_{9}$};   
		\filldraw[black] (.0,1.2) circle (.5pt) node[anchor=south] {$p_{7}$}; 		\filldraw[black] (.4,1.2) circle (.5pt) node[anchor=south]{$p_{8}$}; 				
		\end{tikzpicture}
		\caption{$\beta_{24}=246$}\ 
	\end{subfigure}
	\quad		
	\begin{subfigure}[t]{4.4cm}
		\centering	
		\begin{tikzpicture}[scale=2]
		\draw[] (.4,1.2)--(0,0)--(0,.6)--(0,1.2)--(.8,.6) 	(.4,1.2)--(.4,.6)--(.4,0)	(.8,0)--(.8,.6)--(.8,1.2)--(.4,0)  ;
		\filldraw[black] (0,0) circle (.5pt) node[anchor=north]  {$p_1$}; 
		\filldraw[black] (.8,0) circle (.5pt) node[anchor=north] {$p_3$}; 
		\filldraw[black] (.4,0) circle (.5pt) node[anchor=north] {$p_2$}; 
		\filldraw[black] (0,.6) circle (.5pt) node[anchor=east]  {$p_4$};   
		\filldraw[black] (.4,.6) circle (.5pt) node[anchor=east]  {$p_5$}; 		
		\filldraw[black] (.8,.6) circle (.5pt) node[anchor=west]  {$p_{6}$}; 
		\filldraw[black] (.8,1.2) circle (.5pt) node[anchor=south] {$p_{9}$};   
		\filldraw[black] (.0,1.2) circle (.5pt) node[anchor=south] {$p_{7}$}; 		\filldraw[black] (.4,1.2) circle (.5pt) node[anchor=south]{$p_{8}$}; 				
		\end{tikzpicture}
		\caption{$\beta_{24}=180$}\ 
	\end{subfigure}
	\quad	
	\begin{subfigure}[t]{4.2cm}
		\centering	
		\begin{tikzpicture}[scale=2]
		\draw[] (.4,1.2)--(0,0)--(0,.6)--(0,1.2)--(.8,0) 	(.4,1.2)--(.4,.6)--(.4,0)--(.8,.6)	(.8,0)--(.8,.6)--(.8,1.2)--(.4,.75)  ;
		\filldraw[black] (0,0) circle (.5pt) node[anchor=north]  {$p_1$}; 
		\filldraw[black] (.8,0) circle (.5pt) node[anchor=north] {$p_3$}; 
		\filldraw[black] (.4,0) circle (.5pt) node[anchor=north] {$p_2$}; 
		\filldraw[black] (0,.6) circle (.5pt) node[anchor=east]  {$p_4$};   
		\filldraw[black] (.4,.75) circle (.5pt) node[anchor=west]  {$p_5$}; 		
		\filldraw[black] (.8,.6) circle (.5pt) node[anchor=west]  {$p_{6}$}; 
		\filldraw[black] (.8,1.2) circle (.5pt) node[anchor=south] {$p_{9}$};   
		\filldraw[black] (.0,1.2) circle (.5pt) node[anchor=south] {$p_{7}$}; 		\filldraw[black] (.4,1.2) circle (.5pt) node[anchor=south]{$p_{8}$}; 				
		\end{tikzpicture}
		\caption{$\beta_{24}=237$}\label{fig43}
	\end{subfigure}			
	\caption{}\label{fig36}
\end{figure}

\begin{Discussion}\label{weakenhypo}\normalfont
Here we answer the following question: what happens if we weaken the hypothesis of Theorem~\ref{connected}?  Let $P$ be a poset as defined in  Theorem~\ref{connected}. When $P$ is disconnected , it follows from \cite[Corollary\ 3.2, Theorem\ 3.13]{VEE22} that $\beta_{24}(R[\MI(P)]) = 0$ if and only if $P$ is a disjoint union of two posets $P_1$ and  $P_2$ such that $\MI(P_1)$ is a planar distributive lattice with $\beta_{24}(R[\MI(P_1)]) = 0$ and $P_2$ is a chain.. 

On the other hand, suppose that $P$ is connected and there exists a minimal  element of $P$ which is covered by a maximal element. Using the proof of  Theorem~\ref{connected}, we may replace the poset $P$ by $P_n$ and assume that $n_i=3$ for all $1\leq i \leq 3$. Let $X'$ be as defined in Notation~\ref{construction1}. Observe that $\# X'\in  \{1, 2,3\}$. If  $\# X'=1$ or 2, then $\beta_{24}(R[\MI(P)]) \neq 0$ by the argument of the proof  of Theorem~\ref{connected}.\par 

Now, consider the case when $\# X'= 3$. We know that if $P$ is as shown in figure~\ref{fig30}, then $\beta_{24}(R[\MI(P)]) = 0$. So we may assume that $P$ is not as shown in figure~\ref{fig30}. By Theorem~\ref{planar}, we may assume that there is no subposet $P'$ of $P$, as defined in Theorem~\ref{planar}, with $\beta_{24}(R[\MI(P')]) \neq 0$. By Lemma~\ref{covering3elements}, we may assume that there is no element in $P$ such that either it cover three elements or it is covered by three elements. By Lemma~\ref{subposetconnected}, we may assume that the subposets of $P$ defined on the underlying sets $P\setminus \{p_{1,1},p_{2,1}, p_{3,1}\}$ and $P\setminus \{p_{1,3},p_{2,3}, p_{3,3}\}$ are not connected. Then $P$ is isomorphic to one of the posets as shown in Figure~\ref{fig36}. One can use a computer to check that if $P$ is isomorphic to one of the posets as shown in Figure~\ref{fig40}-\ref{fig41}, then $\beta_{24}(R[\MI(P)]) = 0$ otherwise $\beta_{24}(R[\MI(P)]) \neq0$.
\end{Discussion} 
 
\begin{Remark}\label{reduction3chains}\normalfont
 Let $P$ be a poset. Let $\mathcal{S}=\mathop{\cup}_{i=1}^{3}\{p_{i,1},\ldots ,p_{i,n_i}\}$ be a subset of the underlying set of $P$ such that
 \begin{enumerate}
 	\item[$(a)$] $p_{1,1},p_{2,1},p_{3,1}$ are distinct and $p_{1,n_{1}},p_{2,n_{2}},p_{3,n_{3}}$ are distinct;
 	\item[$(b)$]
 	$B :=\{p_{1,1},p_{2,1}, p_{3,1}\}$ and  $B' := \{p_{1,n_{1}},p_{2,n_2},p_{3,n_3}\}$ are antichains in $P$;
 	\item[$(c)$] for all $1\leq i \leq 3$, $n_i\geq 3$; $\{p_{i,1},\ldots ,p_{i,n_i} \}$  is a chain in $P$ with $p_{i,1}\lessdot\cdots\lessdot p_{i,n_i}$.
 \end{enumerate}
Using Discussion~\ref{reducesubposet} and the arguments of the proof  of Theorem~\ref{planar}, we  can reduce $P$ to the poset $Q^{\partial}$, where $Q^{\partial}$ is a poset on the underlying set $\mathcal{S}$. Note that $B$ and $B'$ are the sets of minimal and maximal elements of $Q^{\partial}$ respectively. by Discussion~\ref{reducesubposet}, $R[\MI(P)]$ does not satisfy property $N_2$ if $R[\MI(Q^\partial)]$ does not satisfy property $N_2$ which can be easily checked using Theorem~\ref{connected} and Discussion~\ref{weakenhypo}.
\end{Remark}

\section{Property $N_p$ of Hibi rings for $p\geq3$}\label{propertn3}

\begin{figure}[h]
	\centering	
	\begin{tikzpicture}[scale=2]
	\draw[](0,1.2)--(.6,0);
	\draw[dotted](0,0)--(0,1.2) (.6,0)--(.6,0.6)--(.6,1.2);
\filldraw[black] (0,0) circle (.5pt) node[anchor=north] {$p_1$};	
\filldraw[black] (0.6,0) circle (.5pt) node[anchor=north] {$q_1$};	
\filldraw[black] (0,1.2) circle (.5pt) node[anchor=south] {$p_n$};
\filldraw[black] (.6,1.2) circle (.5pt) node[anchor=south] {$q_m$};			
	\end{tikzpicture}
	\caption{$P_{n,m}$; $n,m\geq 2$}\label{fig14}
\end{figure}
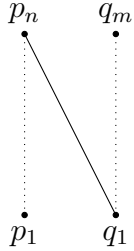	
\begin{Lemma}\label{beta35example}
Let $P_{n,m}$, where $n,m\geq 2$, be the poset as shown in Figure~\ref{fig14}. Then $R[\MI(P_{n,m})]$ does not satisfy $N_3$.
\end{Lemma}

\begin{proof}
Observe that $R[\MI(P_{n,m})]$ satisfies property $N_2$, by Theorem~\ref{ene}. Let $x= p_2$ and $y=p_n$. Reduce $P_{n,m}$ to $P_{2,m}$ using the methods of Notation~\ref{construction1}. Now in $P_{2,m}$, let $x = q_2$ and $y=q_m$. Reduce $P_{2,m}$ to $P_{2,2}$ using the method discussed in Notation~\ref{construction1}. For $n,m=2$, one can use a computer to check that $\beta_{35}(R[\MI(P_{n,m})]) \neq 0$. By Lemma~\ref{reduceposet} and Proposition~\ref{homologicallypure}, we have $\beta_{35}(R[\MI(P_{n,m})]) \neq 0$. This completes the proof.
\end{proof} 
 
\begin{Lemma}\label{beta35planar}
 Let $P$ be a poset such that $\MI(P)$ is a planar distributive lattice. Assume that $P$ has two minimal and maximal elements. If $R([\MI(P)])$ satisfies property $N_3$, then $P$ is a disjoint union of two chains. 
\end{Lemma}
 
\begin{proof}
 Suppose that $R([\MI(P)])$ satisfies property $N_3$. Then, it also satisfies property $N_2$. So $P$ is simple otherwise there exists an element $p\in P$ which is comparable to every element of $P$. By hypothesis, $p$ is neither a minimal element nor a maximal element. Let $P_1=\{q\in P: q<p\}$ and $P_2=\{q\in P: q>p\}$. Since $P_1$ and $P_2$ are not chains, $R([\MI(P)])$ does not satisfy property $N_2$ by Lemma~\ref{nonsimple}, which is a contradiction. By  Corollary~\ref{N2planarposets}, $P$ is isomorphic to one of the posets as shown in Figure~\ref{fig20}. If $P$ is not isomorphic to the poset shown in Figure~\ref{fig15}, then it will contain a cover-preserving subposet as shown in Figure~\ref{fig14}, call  it $P'$. Let $B$ and $B'$ be the sets of minimal  and maximal elements of $P'$ respectively. Hence, by Discussion~\ref{subsemigroup} and Lemma~\ref{beta35example}, $\beta_{35}(R[\MI(P)]) \neq 0$. This concludes the proof.
\end{proof}	

Now we prove our main theorem about property $N_3$ of Hibi rings associated to connected posets.

\begin{Theorem}\label{beta35connected}
  Let $P$ be a connected poset. Assume that $P$ has atleast two minimal and atleast two maximal elements. Then $R[\MI(P)]$ does not satisfy property $N_3$.
\end{Theorem}
 
 \begin{proof}
 $Claim:$ There exist two  maximal chains $C_1=\{p_1,\ldots, p_r\}$ and $C_2=\{q_1,\ldots, q_s\}$ of $P$ such that $p_1\lessdot\cdots\lessdot p_r$, $q_1\lessdot\cdots\lessdot q_s$, $p_1\neq q_1$ , $p_r\neq q_s$ and $r,s\geq2$.\par
 Assume the claim. Let $\mathcal{S}=C_1 \cup C_2$. Using Discussion~\ref{reducesubposet} and the proof  of Theorem~\ref{planar}, we  can reduce $P$ to the poset $Q^{\partial}$, where $Q^{\partial}$ is a poset on the underlying set $\mathcal{S}$ and it is enough to show that $R[\MI(Q^{\partial})]$ does not satisfy property $N_3$. Observe that $Q^{\partial}$ is connected, $\{p_1,q_1\}$ and $\{p_r,q_s\}$ are the sets of minimal and maximal elements of $Q^{\partial}$ respectively. By Lemma~\ref{beta35planar}, $R[\MI(Q^{\partial})]$ does not satisfy property $N_3$. This completes the proof.\par
 
Now we prove the claim. Let $C$ be a maximal chain in $P$ with the minimal element $p$ and maximal element $q$. Fix a maximal element $q'\in P$ where $q'\neq q$. If there exists a maximal chain $C'$ with the maximal element $q'$ and the minimal element not equal to $p$, then we are done. So we may assume that all maximal chains with the maximal element $q'$ have minimal element $p$. Fix a minimal element $p'\in P$ where $p' \neq p$. If there exists a maximal chain $C''$ with the minimal element $p'$ and maximal element not equal to $q$, then we are done. So we may assume that all maximal chains with the minimal element $p'$ have maximal element $q$. Then, we can take $C_1$ to be a maximal chain from $p$ to $q'$ and $C_2$ to be a maximal chain from  $p'$ to $q$. Hence the proof.
\end{proof}

\begin{Definition}
Let $P$ be a poset and $I_{\MI(P)}$ be the Hibi ideal associated to $\MI(P)$. Then $I_{\MI(P)}$ said to have a {\em linear resolution} if $\beta_{ij}(I_{\MI(P)})=0$ for $j \neq i+2$. We say that the Hibi ring $R[\MI(P)]$ has a linear resolution if $I_{\MI(P)}$ has a linear resolution.
\end{Definition}

\begin{figure}[h]
	\centering
			\begin{tikzpicture}	[scale=2]	
			\draw[thick] (0,1) --  (2,1);
		    \draw[thick,dotted] (2,1) --  (4,1);
		    \draw[thick,dotted] (2,0) --  (4,0);	
			\draw[thick] (5,1) --  (4,1);
			\draw[thick] (0,0) --  (2,0);
			\draw[thick] (5,0) --  (4,0);
			\draw[thick] (0,0) --  (0,1);
			\draw[thick] (1,0) --  (1,1);
	        \draw[thick] (4,0) --  (4,1);	
		 	\draw[thick] (2,0) --  (2,1);
		 	\draw[thick] (5,0) --  (5,1); 
		 	
		\filldraw[black] (0,1) circle (.5pt) node[anchor=south] {$a_1$};	
		\filldraw[black] (1,1) circle (.5pt) node[anchor=south] {$a_2$};	
		\filldraw[black] (0,0) circle (.5pt) node[anchor=north] {$b_1$};
		\filldraw[black] (1,0) circle (.5pt) node[anchor=north] {$b_2$};	
		
		\filldraw[black] (2,1) circle (.5pt) node[anchor=south] {$a_3$};	
		\filldraw[black] (4,1) circle (.5pt) node[anchor=south] {$a_{n-1}$};	
		\filldraw[black] (2,0) circle (.5pt) node[anchor=north] {$b_3$};
		\filldraw[black] (4,0) circle (.5pt) node[anchor=north] {$b_{n-1}$};	
		
		\filldraw[black] (5,1) circle (.5pt) node[anchor=south] {$a_{n}$};	
		\filldraw[black] (5,0) circle (.5pt) node[anchor=north] {$b_n$}; 	
		 	   			
			\end{tikzpicture}
	\caption{$L$}\label{fig:linearresolution}
\end{figure}

Recall the notion of graphs from Subsection~\ref{subsec:graphtheory}. The following lemma will be needed in the proof of our main theorem about property $N_p$ for $p\geq 4$.

\begin{Lemma}\label{lem:chordal}
Let $L$ be a distributive lattice as shown in Figure~\ref{fig:linearresolution}. Then the comparability graph $G_L$ of $L$ is chordal. 	
\end{Lemma}
\begin{proof}
	First break the underlying set of $L$ in two disjoint subsets $A_1 = \{a_1,\ldots,a_n\}$ and $A_2 =\{b_1,\ldots,b_n\}$ (see Figure~\ref{fig:linearresolution} for notational conventions). Let $C=(c_1,\ldots,c_r)$ be a induced cycle of $G_L$ of length $\geq 4$. If  $\{c_1,\ldots,c_r\}\cap A_i \geq 3$ for any $i\in \{1,2\}$, then $C$ has a chord because every pair in $A_i$ is an edge of $G_L$. So we may assume that $r =4$ and $\#(\{c_1,\ldots,c_r\}\cap A_i) =2$ for all $i$. Let $\{c_{i_1},c_{i_2}\}\subseteq A_1$ and $\{c_{i_3},c_{i_4}\}\subseteq A_2$. Without loss of generality, we may assume that $c_1= c_{i_1}$ and $c_1<c_{i_2}$ in $L$. Let $c\in \{c_{i_3},c_{i_4}\}$ be such that $\{c_1,c\}$ is an edge in $C$. Therefore, $c_1$ and $c$ are comparable in $L$; therefore $c<c_1$ because $c_1\in A_1$ and $c\in A_2$. Therefore $c<c_{i_2}$. Hence $(c_1,c,c_{i_2})$ is a induced chain in $G_L$. Thus $C$ has a chord. This completes the proof.
\end{proof}

\begin{Example}\label{twoposets}\normalfont
Let $P_1$ be an antichain of cardinality  three and $P_2$ be a poset such that it is a disjoint union of two chains of length 1. By \cite[\S \ 3, Corollary]{[HIBI87]}, $R[\MI(P_i)]$ is a Gorenstein ring for all $i=1,2$. For all $i=1,2$, the Hibi ring $R[\MI(P_i)]$ is Cohen-Macaulay, it is a quotient of a polynomial ring in $\#\MI(P_i)$ variables and the Krull-dimension of $R[\MI(P_i)]$ is $\#P_i+1$. So the Auslander-Buchsbaum formula implies that $\projdim(R[\MI(P_i)])= \#\MI(P_i)-\# P_i -1$ for $i=1,2$. It is easy to see that $\projdim(R[\MI(P_i)]) =4$ for all $i=1,2$. By self-duality of minimal free resolution of Gorenstein rings, we obtain that  $\beta_{4j} (R[\MI(P_i)])\neq 0$ for some $j\geq 6$ and for all $i=1,2$ irrespective of the characteristic of the field $K$. 
\end{Example}

\begin{Theorem}\label{propertyn4}
Let $P$ be a poset  and $p\geq4$. Let $P'= \{p_{i_1},...,p_{i_r}\}$ be the subset of all elements of $P$ which are comparable to every element of $P$. Let $P''$ be the induced subposet of $P$ on the set $P\setminus P'$. Then the following are equivalent:
\begin{enumerate}
	\item[$(a)$] 
	$R[\MI(P)]$ satisfies property $N_p$.
	\item[$(b)$]
	$R[\MI(P)]$ satisfies property $N_4$.
	\item[$(c)$]
	Either $P$ is a chain or $P''$ is a disjoint union of a chain and an isolated element.
	\item[$(d)$]
	Either $R[\MI(P)]$ is a polynomial ring or $K[\MI(P'')]/\ini_<(I_{\MI(P'')})$ has a linear resolution.
	\item[(e)]
	Either $R[\MI(P)]$ is a polynomial ring or it has a linear resolution.
\end{enumerate}
\end{Theorem}

Before going to the proof of the theorem, we remark that not all of the equivalent statements are new. For example, $(c)\iff (e)$ was proved in~\cite[Corollary\ 10]{EQR13regularityplanardistlattices} and $(e)\implies (d)$ follows from \cite[Corollary\ 2.7]{CV20squarefree}.

\begin{proof}
	$(a)\implies(b)$ is trivial.\par
	
	$(b) \implies (c)$  If $\width(P)\geq 3$, then there exists an antichain $P_1$ in $P$ of cardinality three. By Discussion~\ref{subsemigroup}, $\beta_{ij}(R[\MI(P_1)])\leq \beta_{ij}(R[\MI(P)])$ for all $i$ and $j$. Since $\beta_{4j}(R[\MI(P_1)])\neq 0$ for some $j\geq 6$ by Example~\ref{twoposets},  $\beta_{4j}(R[\MI(P)])\neq 0$. Thus, $R[\MI(P)]$ does not satisfy property $N_4$. So we may assume that $\width(P)\leq 2$. If $\width(P)=1$, then $P$ is a chain. We now consider $\width(P)= 2$. Observe that $P''$ is simple. Since $R[\MI(P'')]$ satisfies property $N_4$, it also satisfies property $N_3$. By Lemma~\ref{beta35planar}, $P''$ is a disjoint union of two chains. Suppose that $P''$ is a poset on the set $\mathop{\cup}_{i=1}^{2}\{p_{i,1},\ldots ,p_{i,n_i}\}$ such that $\{p_{i,1},\ldots ,p_{i,n_i} \}$ is a chain in $P''$ with $p_{i,1}\lessdot\cdots\lessdot p_{i,n_i}$ for all $i=1, 2$. We have to show that either $n_1=1$ or $n_2=1$. Suppose that, on the contrary, $n_i\geq2$ for all $i=1,2$. Let $P_2$ be the induced subposet of $P''$ on the set $\mathop{\cup}_{i=1}^{2}\{p_{i,1} ,p_{i,2}\}$. Let $B$ and $B'$ be the sets of minimal and maximal elements of $P_2$ respectively. By Example~\ref{twoposets} and Discussion~\ref{subsemigroup}, $\beta_{4j}(R[\MI(P)]) \neq 0$ for some $j\geq 6$ which is a contradiction. Hence the proof.\par
	
	$(c) \implies (d)$ If $P$ is a chain, then $R[\MI(P)]$ is a polynomial ring. Observe that the distributive lattice $\MI(P'')$ is as shown in Figure~\ref{fig:linearresolution}.  The ideal $\ini_<(I_{\MI(P'')})$ is the Stanley-Reisner ideal of the order complex $\Delta(\MI(P''))$ of $\MI(P'')$ (see \cite[Subsection\ 4.1]{VEE22}). It was observed in Subsection~\ref{subsec:graphtheory} that $\Delta(\MI(P''))=\Delta(G_{\MI(P'')})$ where $G_{\MI(P'')}$ is the comparability graph of ${\MI(P'')}$. Now the result follows from Lemma~\ref{lem:chordal} and \cite[Theorem\ 1]{Fro90}.
	
	$(d) \implies (e)$ Since the Betti numbers of $K[\MI(P'')]/\ini_<(I_{\MI(P'')})$ over the ring $K[\MI(P'')]$ are greater than equal to those of $R[\MI(P'')]$~\cite[Theorem\ 22.9]{[PEEVA10]}, we get that $R[\MI(P'')]$ has a linear resolution. Thus, $R[\MI(P)]$ has a linear resolution by Corollary~\ref{nonsimpleelements}.\par
	
	$(e) \implies(a)$ is immediate.
\end{proof}

We now use Theorem~\ref{propertyn4} and \cite[Theorem\ 1]{Fro90}
to characterize comparability graph of distributive lattices which are chordal. 
It is immediate that for a chain $P$ of length $n$, $G_P$ is the complete graph on the set $[n+1]$ which is chordal. 

\begin{Corollary}\label{coro:chordalcomparabilitygraph}
Let $L=\MI(P)$ be a distributive lattice and $G_L$ be the comparability graph of $L$. For $P$, let $P''$ be as defined in Theorem~\ref{propertyn4}. Then $G_L$ is chordal if and only if $P$ is a chain or $P''$ is a disjoint union of a chain and an isolated element.
\end{Corollary}

\section{Complete intersection Hibi rings}

In this section, we will combinatorially characterize complete intersection Hibi rings.

Let $P_1$ and $P_2$ be two posets and $P$ be the ordinal sum of $P_1$ and $P_2$. Let $R[\MI(P_1)] = K[\{x_\alpha : \alpha \in \MI(P_1) \}]/I_{\MI(P_1)}$, $R[\MI(P_2)] = K[\{y_\beta : \beta \in \MI(P_2) \}]/I_{\MI(P_2)}$ and $R[\MI(P)] = K[\{z_\gamma : \gamma \in \MI(P) \}]/I_{\MI(P)}$.

\begin{Lemma}\label{p_1}
	Let $P_1$, $P_2$ and $P$ be as above. Then 
	$$R[\MI(P)] \cong (R[\MI(P_1)] \tensor_{K} R[\MI(P_2)])/ (x_{P_1}-y_\emptyset).$$
\end{Lemma}

\begin{proof}
	Let $T = K[\{x_\alpha\ : \alpha \in \MI(P_1) \}\cup \{y_\beta\ : \beta \in \MI(P_2) \}]/(x_{P_1}-y_\emptyset)$ and $T' = T/(I_{\MI(P_1)}T  + I_{\MI(P_2)}T)$. Define a map
	$$\varphi : K[\MI(P)]\to T$$ by  
	\[ \varphi(z_\gamma) = \begin{cases} 
	x_\gamma & \text{if} \quad \gamma \subseteq P_1, \\
	y_{\gamma'} & \text{if} \quad \gamma = P_1 \cup \gamma',\ \text{where} \
	\gamma' \subseteq P_2 . \\
	\end{cases}
	\] 
	It is easy to see that $\varphi$ is an isomorphism. If $\alpha, \beta \in \MI(P)$ are incomparable then either $\alpha, \beta \in \MI(P_1)$ or $\alpha = P_1 \cup \alpha'$ and $\beta = P_1 \cup \beta'$ where $\alpha',\beta' \in \MI(P_2)$ and $\alpha',\beta'$ incomparable. Let $\pi : T \to T'$ be the natural projection. Thus, $\pi \circ \varphi :  K[\MI(P)]\to T'$ and $\ker(\pi \circ \varphi) = \varphi^{-1}(I_{\MI(P_1)}T  + I_{\MI(P_2)}T)$. 
	
Thus, it is sufficient to show that  $\varphi (I_{\MI(P)}) = I_{\MI(P_1)}T  + I_{\MI(P_2)}T$. The proof of this is similar to the proof of Lemma~\ref{ordinalsum}.
\end{proof}

\begin{Example}\label{beta23}
	Let $P_1$ and $P_2$ be the posets as shown in Figure~\ref{fig11} and Figure~\ref{fig13} respectively. Then the respective graded Betti table of $R[\MI(P_1)]$ and $R[\MI(P_2)]$ are the following:\par
	
	$\begin{matrix}
	&0&1&2&3&4\\\text{total:}&1&9&16&9&1\\\text{0:}&1&\text{.}&\text{.}&\text{.}&\text{.}\\\text{1:}&\text{.}&9&16&9&\text{.}\\\text{2:}&\text
	{.}&\text{.}&\text{.}&\text{.}&1\\\end{matrix}$ \quad
	\quad
	$\begin{matrix}
	&0&1&2\\\text{total:}&1&3&2\\\text{0:}&1&\text{.}&\text{.}\\\text{1:}&\text{.}&3&2\\\end{matrix}
	$
\end{Example}

\begin{figure}[h]
	\begin{subfigure}[t]{4cm}
		\centering	
		\begin{tikzpicture}[scale=2] 
		\draw[] (0,0)--(0,.6)--(.7,0)--(.7,.6)--(0,0)
		(0,1.8)--(0,2.4)--(.7,1.8)--(.7,2.4)--(0,1.8);
		\draw[dotted](0,.6)--(0,1.2)--(.7,.6)--(.7,1.2)--(0,.6)
		(.7,1.2)--(.7,1.8)--(0,1.2)--(0,1.8)--(.7,1.2);	
		\filldraw[black] (.7,1.2) circle (.5pt) node[anchor=north]  {}; 
		\filldraw[black] (0,1.2) circle (.5pt) node[anchor=east] {}; 
		
		\filldraw[black] (0,0) circle (.5pt) node[anchor=north]  {$p_1$}; 
		\filldraw[black] (0,.6) circle (.5pt) node[anchor=east] {$p_3$}; 
		\filldraw[black] (.7,0) circle (.5pt) node[anchor=north] {$p_2$}; 
		\filldraw[black] (.7,.6) circle (.5pt) node[anchor=west] {$p_4$}; 
		\filldraw[black] (0,1.8) circle (.5pt) node[anchor=east] {$p_{2n-3}$}; 
		\filldraw[black] (.7,1.8) circle (.5pt) node[anchor=west] {$p_{2n-2}$}; 
		\filldraw[black] (0,2.4) circle (.5pt) node[anchor=south] {$p_{2n-1}$}; 
		\filldraw[black] (.7,2.4) circle (.5pt) node[anchor=south] {$p_{2n}$}; 
		\end{tikzpicture}
		\caption{}\label{fig12}
	\end{subfigure}	 
	\quad	 
	\begin{subfigure}[t]{4cm}
		\centering	
		\begin{tikzpicture}[scale=1]
		\draw[fill= white] ;
		\filldraw[black] (0,0) circle (.8pt) node[anchor=north] {}; 
		\filldraw[black] (.7,0) circle (.8pt) node[anchor=north] {}; 
		\filldraw[black] (1.4,0) circle (.8pt) node[anchor=north]{}; 			
		\end{tikzpicture}
		\caption{$\beta_{23}=16$}\label{fig11}
	\end{subfigure}
	\quad
	\begin{subfigure}[t]{4cm}
		\centering	
		\begin{tikzpicture}[scale=2]
		\draw[] (0,0)--(0,.6);
		
		\filldraw[black] (.5,0) circle (.5pt) node[anchor=north] {$p_2$}; 
		\filldraw[black] (0,0) circle (.5pt) node[anchor=north]  {$p_1$}; 
		\filldraw[black] (0,.6) circle (.5pt) node[anchor=south] {$p_3$};  
		\end{tikzpicture}
		\caption{$\beta_{23}=2$}\label{fig13}
	\end{subfigure}
	\caption{}
\end{figure}

Now, we prove the main theorem of this section.

\begin{Theorem}\label{completeinter}
	Let $P$ be a poset and $P'= \{p_{i_1},...,p_{i_r}\}$ be the subset of all elements of $P$ which are comparable to every element of $P$. Let $P''$ be the induced subposet of $P$ on the set $P\setminus P'$. Then the following are equivalent:\par
	$(a)$ $R[\MI(P)]$ is a complete intersection.\par
	$(b)$ Either $P$ is a chain or $P''$ is as shown in Figure~\ref{fig12}.
\end{Theorem}

\begin{proof}
	$(a)\implies(b).$ Suppose that $R[\MI(P)]$ is a complete intersection. Therefore, $\beta_{23}(R[\MI(P)])= 0$. We break the proof by width of the poset. If $\width(P)\geq 3$, then there exists an antichain $P_1$ of $P$ of cardinality 3. By Discussion~\ref{subsemigroup}, $\beta_{23}(R[\MI(P_1)])\leq \beta_{23}(R[\MI(P)])$. Since $\beta_{23}(R[\MI(P_1)])\neq 0$ by Example~\ref{beta23}, we obtain that $\beta_{23}(R[\MI(P)])\neq 0$. So we may assume that $\width(P)\leq 2$. If $\width(P)=1$, then $P$ is a chain. Hence, the only case we need to consider is $\width(P)= 2$. Now if $P''$ is not as shown in Figure~\ref{fig12}, then $P''$ contain the poset as shown in Figure~\ref{fig13} as a cover-preserving subposet, call it $P_2$. Let $B$ and $B'$ be the sets of minimal and maximal elements of $P_2$ respectively. From Discussion~\ref{subsemigroup} and Example~\ref{beta23}, $\beta_{23}(R[\MI(P)])\neq 0$. This concludes the proof.\par 
	
	$(b)\implies(a).$ If $P$ is a chain, then $R[\MI(P)]$ is a polynomial ring. So we may assume that $P$ is not a chain.  Since $R[\MI(P)] \cong R[\MI(P'')] \tensor_{K}K[y_1,...,y_r]$, it is enough to show that $R[\MI(P'')]$ is a complete intersection. For $1\leq i \leq n$, let $P_i =\{p_{2i-1},p_{2i}\}$ and $Q_i = \{a\in P'': a\leq p_{2i-1}\}\cup\{p_{2i}\}$ be the subposets of $P$. Observe that $Q_n = P''$. For $1\leq i \leq n-1$, by Lemma~\ref{p_1}, 
	$$R[\MI(Q_{i+1})] \cong (R[\MI(Q_i)] \tensor_{K} R[\MI(P_{i+1})])/ (x_{Q_i}-y_\emptyset)$$
	where $\emptyset$ is the minimal element $\MI(P_{i+1})$. 
	We prove the theorem by induction on $i$. It is easy to see that the result holds for $i=1$. Now assume that the result holds for $i$. Since $ R[\MI(P_{i+1})]\cong R[\MI(Q_1)]$, we get $R[\MI(Q_i)] \tensor_{K} R[\MI(P_{i+1})]$ is a complete intersection. Hence, $R[\MI(Q_{i+1})]$ is a complete intersection. Hence the proof.
\end{proof}


\begin{thebibliography}{EHSM15}

\bibitem[ACI15]{[ACI15]}
Luchezar~L. Avramov, Aldo Conca, and Srikanth~B. Iyengar.
\newblock Subadditivity of syzygies of {K}oszul algebras.
\newblock {\em Math. Ann.}, 361(1-2):511--534, 2015.

\bibitem[BH97]{BH}
Winfried Bruns and J\"{u}rgen Herzog.
\newblock Semigroup rings and simplicial complexes.
\newblock {\em J. Pure Appl. Algebra}, 122(3):185--208, 1997.

\bibitem[CV20]{CV20squarefree}
Aldo Conca and Matteo Varbaro.
\newblock Square-free {G}r\"{o}bner degenerations.
\newblock {\em Invent. Math.}, 221(3):713--730, 2020.

\bibitem[DM17]{[DM17]}
Priya Das and Himadri Mukherjee.
\newblock First syzygy of hibi rings associated with planar distributive
  lattices.
\newblock {\em arXiv preprint arXiv:1704.08286}, 2017.

\bibitem[EHH15]{[EHH15]}
Viviana Ene, J\"{u}rgen Herzog, and Takayuki Hibi.
\newblock Linearly related polyominoes.
\newblock {\em J. Algebraic Combin.}, 41(4):949--968, 2015.

\bibitem[EHSM15]{[EHM15]}
Viviana Ene, J\"{u}rgen Herzog, and Sara Saeedi~Madani.
\newblock A note on the regularity of {H}ibi rings.
\newblock {\em Manuscripta Math.}, 148(3-4):501--506, 2015.

\bibitem[Ene15]{[ENE15]}
Viviana Ene.
\newblock Syzygies of {H}ibi rings.
\newblock {\em Acta Math. Vietnam.}, 40(3):403--446, 2015.

\bibitem[EQR13]{EQR13regularityplanardistlattices}
Viviana Ene, Ayesha~Asloob Qureshi, and Asia Rauf.
\newblock Regularity of join-meet ideals of distributive lattices.
\newblock {\em Electron. J. Combin.}, 20(3):Paper 20, 8, 2013.

\bibitem[Fr\"{o}90]{Fro90}
Ralf Fr\"{o}berg.
\newblock On {S}tanley-{R}eisner rings.
\newblock In {\em Topics in algebra, {P}art 2 ({W}arsaw, 1988)}, volume~26 of
  {\em Banach Center Publ.}, pages 57--70. PWN, Warsaw, 1990.

\bibitem[GL86]{[GL86]}
Mark Green and Robert Lazarsfeld.
\newblock On the projective normality of complete linear series on an algebraic
  curve.
\newblock {\em Invent. Math.}, 83(1):73--90, 1986.

\bibitem[HHO18]{[HHO18]}
J\"{u}rgen Herzog, Takayuki Hibi, and Hidefumi Ohsugi.
\newblock {\em Binomial ideals}, volume 279 of {\em Graduate Texts in
  Mathematics}.
\newblock Springer, Cham, 2018.

\bibitem[Hib87]{[HIBI87]}
Takayuki Hibi.
\newblock Distributive lattices, affine semigroup rings and algebras with
  straightening laws.
\newblock In {\em Commutative algebra and combinatorics ({K}yoto, 1985)},
  volume~11 of {\em Adv. Stud. Pure Math.}, pages 93--109. North-Holland,
  Amsterdam, 1987.

\bibitem[HO17]{OH17chordalcomparabilitygraphs}
Takayuki Hibi and Hidefumi Ohsugi.
\newblock A {Gr{\"o}bner} basis characterization for chordal comparability
  graphs.
\newblock {\em Eur. J. Comb.}, 59:122--128, 2017.

\bibitem[Kem90]{[KEMPH90]}
George~R. Kempf.
\newblock Some wonderful rings in algebraic geometry.
\newblock {\em J. Algebra}, 134(1):222--224, 1990.

\bibitem[M2]{[M2]}
Daniel~R. Grayson and Michael~E. Stillman.
\newblock Macaulay2, a software system for research in algebraic geometry.
\newblock Available at \url{http://www.math.uiuc.edu/Macaulay2/}.

\bibitem[Pee11]{[PEEVA10]}
Irena Peeva.
\newblock {\em Graded syzygies}, volume~14 of {\em Algebra and Applications}.
\newblock Springer-Verlag London, Ltd., London, 2011.

\bibitem[Sage]{[sagemath]}
{The Sage Developers}.
\newblock {\em {S}ageMath, the {S}age {M}athematics {S}oftware {S}ystem
  ({V}ersion 9.2)}, 2020.
\newblock {\tt https://www.sagemath.org}.

\bibitem[Sta12]{[STAN12]}
Richard~P. Stanley.
\newblock {\em Enumerative combinatorics. {V}olume 1}, volume~49 of {\em
  Cambridge Studies in Advanced Mathematics}.
\newblock Cambridge University Press, Cambridge, second edition, 2012.

\bibitem[Stu96]{[STURM96]}
Bernd Sturmfels.
\newblock {\em Grobner bases and convex polytopes}, volume~8.
\newblock American Mathematical Soc., 1996.

\bibitem[Vee22]{VEE22}
Dharm Veer.
\newblock {G}reen-{L}azarsfeld property ${N_p}$ for {S}egre product of {H}ibi
  rings.
\newblock {\em Submitted}, 2022.

\end{thebibliography}

\end{document}